\documentclass[11pt,a4paper]{article}
\pdfoutput=1
\usepackage{amssymb}
\usepackage{amsmath}
\usepackage{amsthm}
\usepackage[T1]{fontenc}
\usepackage{textcomp}
\usepackage{times}
\usepackage{graphicx}
\usepackage[all]{xy}
\usepackage[hyperref, colorlinks=true, pdfpagemode=none, pdfmenubar=false, pdfstartview=FitH, linkcolor=blue, citecolor=blue, urlcolor=blue, pdffitwindow=false]{hyperref}
\DeclareMathOperator{\codim}{codim}

\DeclareMathOperator{\img}{Im}

\DeclareMathOperator{\homm}{Hom}
\DeclareMathOperator{\shhomm}{{\it Hom\rm}}
\DeclareMathOperator{\eend}{End}

\DeclareMathOperator{\gl}{GL}
\DeclareMathOperator{\gsl}{SL}
\DeclareMathOperator{\gsp}{Sp}

\DeclareMathOperator{\picr}{Pic}

\DeclareMathOperator{\cl}{Cl}

\DeclareMathOperator{\spec}{Spec\;}

\DeclareMathOperator{\hilb}{Hilb}

\DeclareMathOperator{\ext}{Ext}

\DeclareMathOperator{\tor}{Tor}
\DeclareMathOperator{\supp}{Supp}

\DeclareMathOperator{\ann}{Ann}

\DeclareMathOperator{\qcohcat}{QCoh}
\DeclareMathOperator{\cohcat}{Coh}
\DeclareMathOperator{\modd}{\bf Mod}
\DeclareMathOperator{\lder}{\bf L}
\DeclareMathOperator{\rder}{\bf R}
\DeclareMathOperator{\rank}{rk}
\DeclareMathOperator{\ptfield}{\bf k}

\DeclareMathOperator{\prjdim}{proj\;dim}
\DeclareMathOperator{\depth}{depth}
\DeclareMathOperator{\cxampl}{coh-amp}
\DeclareMathOperator{\torampl}{Tor-amp}

\DeclareMathOperator{\nonorth}{N}

\DeclareMathOperator{\proobjet}{pro}

\DeclareMathOperator{\id}{Id}
\begin{document}
\def\bv{\mathbf{v}}
\def\kgc_{K^*_G(\mathbb{C}^n)}
\def\kgchi_{K^*_\chi(\mathbb{C}^n)}
\def\kgcf_{K_G(\mathbb{C}^n)}
\def\kgchif_{K_\chi(\mathbb{C}^n)}
\def\gpic_{G\text{-}\picr}
\def\gcl_{G\text{-}\cl}
\def\trch_{{\chi_{0}}}
\def\regring{{R}}
\def\regrep{{V_{\text{reg}}}}
\def\givrep{{V_{\text{giv}}}}
\def\lbar{{(\mathbb{Z}^n)^\vee}}
\def\genpx_{{p_X}}
\def\genpy_{{p_Y}}
\def\genpcn_{p_{\mathbb{C}^n}}
\def\gnat{gnat}
\def\twalg{{\regring \rtimes G}}

\theoremstyle{definition}
\newtheorem{definition}{Definition}
\newtheorem*{definition*}{Definition}
\theoremstyle{plain}
\newtheorem{theorem}{Theorem}
\newtheorem*{theorem*}{Theorem}
\newtheorem{conjecture}{Conjecture}
\newtheorem{proposition}{Proposition}
\newtheorem*{proposition*}{Proposition}
\newtheorem{corollary}{Corollary}
\newtheorem*{corollary*}{Corollary}
\newtheorem{lemma}{Lemma}
\newtheorem*{claim*}{Claim}

\numberwithin{equation}{section}

\title{Derived McKay correspondence via pure-sheaf transforms}
\author{Timothy Logvinenko}	

\maketitle

\begin{abstract}
In most cases where it has been shown to exist the derived McKay 
correspondence $D(Y) \xrightarrow{\sim} D^G(\mathbb{C}^n)$ can
be written as a Fourier-Mukai transform which sends point sheaves
of the crepant resolution $Y$ to pure sheaves in $D^G(\mathbb{C}^n)$.
We give a sufficient condition for 
$E \in D^G(Y \times \mathbb{C}^n)$ to be the defining object of
such a transform. We use it to construct the first example of
the derived McKay correspondence for a non-projective 
crepant resolution of $\mathbb{C}^3/G$. Along the way we extract
more geometrical meaning out of the Intersection Theorem and
learn to compute $\theta$-stable families of $G$-constellations
and their direct transforms.
\end{abstract}

\section{Introduction} \label{section-intro}

It was observed by McKay in \cite{McKay-GraphsSingularitiesAndFiniteGroups}
that the representation graph (better known now as the \it McKay quiver\rm) 
of a finite subgroup $G$ of $\gsl_2(\mathbb{C})$ is the Coxeter graph 
of one of the affine Lie algebras of type ADE, 
while the configuration of irreducible exceptional divisors
on the minimal resolution $Y$ of $\mathbb{C}^2/G$ is dual to
the Coxeter graph of the finite-dimensional Lie algebra of the same
type. It followed that the subgraph of nontrivial irreducible 
representations coincided with the graph of irreducible exceptional divisors.
This led Gonzales-Sprinberg and Verdier in 
\cite{GsV-ConstructionGeometriqueDeLaCorrespondanceDeMcKay}
to construct an isomorphism of the $G$-equivariant $K$-theory of $\mathbb{C}^2$ 
to the $K$-theory of $Y$, which induced naturally 
a choice of such bijection. This became known as 
\it the (classical) McKay correspondence\rm.

In \cite{Kinosaki-97} M.Reid proposed that the $K$-theory isomorphism 
might lift to the level of derived categories. It became known as 
\it the derived McKay correspondence \rm conjecture: 
\begin{conjecture} \label{conj-derived-mckay-correspondence}
Let $G$ be a finite subgroup of $\gsl_n(\mathbb{C})$ and let $Y$
be a crepant resolution of $\mathbb{C}^n/G$, if one exists. 
Then 
\begin{align} \label{eqn-derived-mckay-correspondence}
D(Y) \xrightarrow{\sim} D^G(\mathbb{C}^n)
\end{align}
where $D(Y)$ and $D^G(\mathbb{C}^n)$ are bounded derived categories of 
coherent sheaves on $Y$ and of $G$-equivariant coherent sheaves on
$\mathbb{C}^n$, respectively.
\end{conjecture}
To date and to the extent of our knowledge this conjecture
has been settled for the following situations: 
\begin{enumerate}
\item \label{case-mckay-bkr} $G \subset \gsl_{2,3}(\mathbb{C})$; $Y$ the
distinguished crepant resolution $G$-$\hilb$; \\
(\cite{KapranovVasserot-KleinianSingularitiesDerivedCategoriesAndHallAlgebras}, Theorem
1.4; \cite{BKR01}, Theorem 1.1).  
\item \label{case-mckay-crawishii} $G \subset \gsl_{3}(\mathbb{C})$ abelian; $Y$ any projective
crepant resolution; \\ (\cite{Craw-Ishii-02}, Theorem 1.1). 
\item \label{case-mckay-kawamata} $G \subset \gsl_{n}(\mathbb{C})$ abelian; 
$Y$ any projective crepant resolution; \\
(\cite{Kawamata-LogCrepantBirationalMapsAndDerivedCategories}, special
case of Theorem 4.2).
\item \label{case-mckay-kaledin} $G \subset \gsp_{2n}(\mathbb{C})$; $Y$ any symplectic (crepant) 
resolution; \\ (\cite{Kal-Bezr-04}, Theorem 1.1).
\end{enumerate}
In the case \ref{case-mckay-kawamata} the construction is not direct
and it isn't clear what form does the equivalence
\eqref{eqn-derived-mckay-correspondence} take, but in each of the
cases \ref{case-mckay-bkr},
\ref{case-mckay-crawishii} and \ref{case-mckay-kaledin}, 
the equivalence \eqref{eqn-derived-mckay-correspondence} is constructed
directly and we observe that the constructed functor sends point 
sheaves $\mathcal{O}_y$ of $Y$ to pure sheaves 
(i.e. complexes with cohomologies concentrated in degree zero) 
in $D^G(\mathbb{C}^n)$. Another property (cf. though
\cite{Orlov-EquivalencesOfDerivedCategoriesAndK3Surfaces}, Theorem
2.18) that these functors share is that each can be written as 
a Fourier-Mukai transform $\Phi_{E}(- \otimes \rho_0)$ (see Def. 
\ref{defn-integral-transform}) for some object
$E \in D^G(Y \times \mathbb{C}^n)$. 

A straightforward application (Prop. \ref{prps-psps-necessary-criterion}) 
of the established machinery of Fourier-Mukai transforms shows that if 
an equivalence $\eqref{eqn-derived-mckay-correspondence}$ is 
a Fourier-Mukai transform $\Phi_{E}(- \otimes \rho_0)$ which
sends point sheaves to pure sheaves, then its defining object $E$
is itself a pure sheaf. Moreover, the fibers of $E$ over $Y$ have 
to be simple ($G$-$\eend_{\mathbb{C}^n}(E_{|y}) = \mathbb{C}$ for all 
$y \in Y$), orthogonal in all degrees ($G$-$\ext^i_{\mathbb{C}^n}(E_{|y_1},
E_{|y_2}) = 0$ if $y_1 \neq y_2$) and the Kodaira-Spencer maps have
to be isomorphisms. 

Let $Y$ now be any irreducible separated scheme of finite type 
over $\mathbb{C}$. A $\gnat$-family $\mathcal{F}$ on $Y$ is a coherent
$G$-sheaf on $Y \times \mathbb{C}^n$, flat over $Y$, such 
that for any $y \in Y$ the fiber $\mathcal{F}_{|y}$ of $\mathcal{F}$ is 
a \it $G$-constellation \rm supported on a single $G$-orbit.
That is, $\mathcal{F}_{|y}$ is a finite length coherent $G$-sheaf on 
$\mathbb{C}^n$ whose support is a single $G$-orbit and whose 
global sections have $G$-representation structure of the regular
representation. Such family $\mathcal{F}$ has a well-defined 
Hilbert-Chow morphism $\pi_\mathcal{F}: Y \rightarrow \mathbb{C}^n/G$,
it sends any $y \in Y$ to the $G$-orbit that $\mathcal{F}_{|y}$ is supported 
on (Prop. \ref{prps-hilbert-chow-map-is-a-morphism}). Let 
$Y$ and $\mathcal{F}$ be any such for which $\pi_\mathcal{F}$ is birational 
and proper. In this paper we give a sufficient condition 
for the functor $\Phi_{\mathcal{F}}(- \otimes \rho_0)$ to be an 
equivalence $\eqref{eqn-derived-mckay-correspondence}$.
Notable, in the view of Prop. \ref{prps-psps-necessary-criterion}, 
is that this condition only asks for the non-orthogonality locus 
of $\mathcal{F}$ to be of high enough codimension. The simplicity 
of $\mathcal{F}$ and the Kodaira-Spencer maps being isomorphisms 
follow automatically:

\begin{theorem} \label{theorem-main-theorem}
Let $G$ be a finite subgroup of $\gsl_n(\mathbb{C})$. Let $Y$ be an
irreducible separated scheme of finite type over $\mathbb{C}$ 
and $\mathcal{F}$ be a $\gnat$-family on $Y$. Assume $Y$ and
$\mathcal{F}$ such that the Hilbert-Chow morphism $\pi_\mathcal{F}$ 
is birational and proper. 

If for every $0 \leq k < (n+1)/2$, the codimension of the subset 
\begin{align} \label{eqn-codimension-of-nonorthogonality}
\nonorth_k = \overline{\{ (y_1,y_2) \in Y \times Y \setminus \Delta  \;\;|\;\; 
G\text{-}\ext_{\mathbb{C}^n}^k
(\mathcal{F}_{|y_1},\mathcal{F}_{|y_2}) \neq 0\}}
\end{align}
in $Y \times Y$ is at least $n + 1 - 2k$, then the functor 
$\Phi_{\mathcal{F}}(- \otimes \rho_0)$ is an equivalence 
of categories $D(Y) \xrightarrow{\sim} D^G(\mathbb{C}^n)$.
\end{theorem}

Once $\Phi_{\mathcal{F}}(- \otimes \rho_0)$ is known to be an equivalence 
usual methods 
(\cite{Roberts-MultiplicitiesAndChernClassesInLocalAlgebra}, Theorem 6.2.2
 and \cite{BKR01}, Lemma 3.1) apply to show that $Y$ is non-singular 
and $\pi_\mathcal{F}$ is crepant. The set $N_k$ in
\eqref{eqn-codimension-of-nonorthogonality} can be thought of as the
locus of the degree $k$ non-orthogonality in $\mathcal{F}$. 

Our proof of Theorem \ref{theorem-main-theorem} is based on 
the ideas introduced 
in \cite{BonOr95} and \cite{BKR01},
particularly on the Intersection Theorem trick introduced in the
latter. However, not wishing to restrict ourselves to just
quasi-projective schemes necessitates more work in applying the
Intersection Theorem. This is done in Section 
\ref{section-cohomological-width}, which is a self-contained piece
of abstract derived category theory for a locally noetherian scheme 
$X$. There we propose a generalisation 
of the concept of the \it homological dimension \rm of 
$E \in D^b_{\text{coh}}(X)$
which we call \it $\tor$-amplitude \rm, and use it to show that the inequality 
$$ \text{ hom. dim. } E \geq \codim_X \supp E $$ of 
\cite{BridgelandMacioca-FourierMukaiTransformsFor_K3_AndEllipticFibrations},
Corollary 5.5 refines to
\begin{align*}
\torampl E \geq \codim_X \supp E + \cxampl E. 
\end{align*}
 
Other notable points of our proof of Theorem
\ref{theorem-main-theorem} are a different approach to Grothendieck 
duality when constructing the left adjoint to
$\Phi_\mathcal{F}(- \otimes \rho_0)$ and
an application of 
$\cite{Logvinenko-Natural-G-Constellation-Families}$, Prop. 1.5
which states that outside the exceptional set of $Y$ any
$\gnat$-family has to be locally isomorphic to the universal family of
$G$-clusters. The latter is everywhere simple and its Kodaira-Spencer
maps are isomorphisms. Then the locus of points of $Y$ where 
objects of $\mathcal{F}$ are not simple or the Kodaira-Spencer map
isn't an isomorphism turns out to have too high a codimension 
to exist at all.

The question of an existence of a derived McKay correspondence which 
sends point sheaves to pure sheaves is thus reduced to that of
an existence of a $\gnat$-family satisfying the non-orthogonality
condition of Theorem $\ref{theorem-main-theorem}$. This is
particularly relevant whenever $G$ is abelian, for then all 
the $gnat$-families on a given resolution $Y \rightarrow \mathbb{C}^n/G$ 
had been classified and their number was shown to be finite and non-zero 
(\cite{Logvinenko-Natural-G-Constellation-Families}, Theorem 4.1). 

When $n = 3$, Theorem $\ref{theorem-main-theorem}$ reduces to:

\begin{corollary} \label{cor-3d-application-of-main-theorem}
Let $G$ be a finite subgroup of $\gsl_3(\mathbb{C})$. Let $Y$, 
$\mathcal{F}$ and $\pi_\mathcal{F}$ be as in Theorem 
$\ref{theorem-main-theorem}$. Let $E_1, \dots, E_k$ be
the irreducible exceptional surfaces of $\pi_\mathcal{F}$.  
Then if general points of any surface $E_i$ are
orthogonal in degree $0$ in $\mathcal{F}$ to general points of any
surface $E_j$ (including case $j = i$) and of any curve 
$E_l \cap E_m$, then $\Phi_\mathcal{F}(- \otimes \rho_0)$ is 
an equivalence of categories.
\end{corollary}
By a general point of an intersection of $k$ exceptional surfaces we 
mean a point that doesn't lie on an intersection of any $k+1$ 
exceptional surfaces. 

In Section $\ref{section-degree-zero-orthogonality}$ we show how to 
compute the degree $0$ non-orthogonality locus of a $\gnat$-family. 
We use this in Section \ref{section-non-projective-example}
to give following application of Corollary 
\ref{cor-3d-application-of-main-theorem}: for $G$ the abelian 
subgroup of $\gsl_3(\mathbb{C})$ known as $\frac{1}{6}(1,1,4) \oplus
\frac{1}{2}(1,0,1)$ (see Section \ref{subsection-the-group}) and for
$Y$ a certain non-projective crepant resolution 
of $\mathbb{C}^3/G$ (see Section \ref{subsection-toric-generalities})
we construct a $\gnat$-family $\mathcal{F}$ on $Y$ which satisfies the condition
in Corollary $\ref{cor-3d-application-of-main-theorem}$. This gives
the first example of the derived McKay correspondence
for a non-projective crepant resolution of $\mathbb{C}^3/G$.

It also leads to an important observation: the properties that 
$\mathcal{F}$ must then possess in view of Proposition 
\ref{prps-psps-necessary-criterion} imply that $Y$ is 
a fine moduli space of $G$-constellations, representing 
the functor of all $\gnat$-families whose members (fibres over
closed points) are isomorphic to members of $\mathcal{F}$. 
At present the only moduli functors known for $G$-constellations 
come from the notion of $\theta$-stability. Their fine moduli spaces 
$M$ (cf. \cite{Craw-Ishii-02}) are constructed via 
the method introduced by King in \cite{King94}. However, $Y$ can't 
be one of $M$
as these are all, due to the GIT nature of their construction in
\cite{King94}, projective over $\mathbb{C}^n/G$. This raises
the question as to whether there could exist a more general notion
of `stability', related perhaps to Bridgeland-Douglas stability
\cite{Bridgeland-StabilityConditionsOnTriangulatedCategories},
which would allow for functors with non-projective moduli spaces.

\bf Acknowledgements: \rm The author would like to express his 
gratitude to S. Mukai, D. Kaledin, D. Orlov and A. Bondal 
for useful discussions while the paper was written, to A. Craw 
for observing a crucial link with the work in 
\cite{Craw-Maclagan-Thomas-05-I}, \cite{Craw-Maclagan-Thomas-05-II} 
which inspired the Proposition \ref{prps-stability-criterion} and to
A. King and an anonymous referee for many helpful comments on
the first draft. The paper was originally completed during the
author's stay at RIMS, Kyoto, and he would like to thank everyone at
the institute for their hospitality. A substantial revision was then
carried out during the author's stay at Mittag-Leffler Institute/KTH,
Stockholm, and he would like to thank them also. 

\section{Cohomological and Tor amplitudes} 
\label{section-cohomological-width}

We clarify terminology and introduce notation. By a point 
of a scheme we mean both a closed and non-closed point unless 
specifically mentioned otherwise. Given a point $x$ on a scheme $X$
we write $(\mathcal{O}_x,\mathfrak{m}_x)$ for the local ring 
of $x$, $\ptfield(x)$ for the residue field $\mathcal{O}_x/\mathfrak{m}_x$ 
and $\iota_x$ for the point-scheme inclusion 
$\spec \ptfield(x) \hookrightarrow X$. Given an irreducible closed
set $C \subset X$, we write $x_C$ for the generic point of $C$ and
we sometimes write simply 
$(\mathcal{O}_C, \mathfrak{m}_C)$ for the local ring of $x_C$. 
All complexes are cochain complexes. Given a right (resp. left) exact
functor $F$ between two abelian categories $\mathcal{A}$ and
$\mathcal{B}$, we denote by $\lder F$ (resp. $\rder F$) the left (resp.
right) derived functor between the appropriate derived categories, 
if it exists, and by $\lder^i F(\bullet)$ (resp. $\rder^i F (\bullet)$)
the $-i$-th cohomology of $\lder F(\bullet)$ (resp. 
the $i$-th cohomology of $\rder F(\bullet)$).

For $X$ a smooth variety the results of Lemmas
\ref{lemma-cohdim-at-least-support-codim} and
\ref{lemma-codim-of-cohdim-p-is-p} below have appeared 
in the proof of Proposition 1.5 in \cite{BonOr95}. We show them 
to hold in a more general setting of a locally noetherian scheme.

\begin{lemma} \label{lemma-cohdim-at-least-support-codim}
Let $X$ be a locally noetherian scheme. Let $\mathcal{F}$ be 
a coherent sheaf on $X$ and $C$ be an irreducible component 
of $\supp_X \mathcal{F}$. Then for every point $x \in C$
\begin{align}
\lder^i \iota_x^* \mathcal{F} \neq 0 \quad \text{ for } 
0 \leq i \leq \codim_X(C). 
\end{align}

\end{lemma}
\begin{proof}
Recall (cf. \cite{Mats86}, \S19) that if a minimal free resolution 
$L_\bullet$ of a finitely generated module $M$ for a local ring
$(R,\mathfrak{m}, k)$ exists, then 
$$\dim_k \tor^i(M, k) = \rank L_i$$ 
Since $X$ is locally noetherian minimal free resolutions of
$\mathcal{F}$ exist in all local rings. Write $F_C$ for 
the localisation of $\mathcal{F}$ to the local ring $\mathcal{O}_C$
of $x_C$. As $\lder^i \iota_x^* \mathcal{F} = \tor^i_{\mathcal{O}_C}(F_C,\ptfield(x))$
it suffices to prove that the length of the minimal 
free resolution of $F_C$ is at least $\codim_X(C)$. 

 Consider 
the standard filtration 
(\cite{Serre-LocalAlgebra}, I, \S7, Theorem 1) of $F_C$ by 
submodules $0 = M_0 \subset \dots \subset M_{n} = F_C$ with 
each $M_i/M_{i-1}$ isomorphic to $\mathcal{O}_C / \mathfrak{p}$ 
for some $\mathfrak{p} \in \supp_{\mathcal{O}_C}(F_C)$. 
As the defining ideal of $C$ is minimal in 
$\supp_X(\mathcal{F})$, $\supp_{\mathcal{O}_C} (F_C)$ consists 
of just $\mathfrak{m}_C$. So each $M_i/M_{i-1}$ is 
isomorphic to $k_C$ and hence $F$ is a finite-length 
$\mathcal{O}_C$-module. Then by the New Intersection Theorem (e.g. \cite{Roberts-MultiplicitiesAndChernClassesInLocalAlgebra},
Theorem 6.2.2) the length of the minimal resolution of $F_C$ is at
least $\dim \mathcal{O}_C$.  As $\dim \mathcal{O}_C= \codim_X(C)$ the
claim follows.
\end{proof}

\begin{lemma} \label{lemma-codim-of-cohdim-p-is-p}
Let $X$ be a locally noetherian scheme. Let $\mathcal{F}$ be a
coherent sheaf on $X$ of finite $\tor$-dimension. For any $p \in
\mathbb{Z}$ define 
\begin{align}
  D_p = \{ x \in X \;|\; \lder^i \iota^*_x \mathcal{F} \neq 0 
\text{ for some } i \geq p \}.
\end{align}

Then each $D_p$ is closed and $\codim_X(D_p) \geq p$. 
\end{lemma}

\begin{proof}
It suffices to prove both claims for the case $X = \spec R$ with 
$R$ noetherian. Write $F$ for $\Gamma(\mathcal{F})$. As 
$\lder^p \iota^*_x \mathcal{F} = \tor_R^{p}(F,\ptfield(x))$
the first claim follows from the upper semicontinuity theorem
(\cite{Grothendieck-EGA-III-2}, \it Th{\'e}or{\`e}me \bf 7.6.9\rm).

For the second claim let $C$ be any irreducible component of $D_p$
and let $F_C$ be the localisation of $F$ to the local ring $\mathcal{O}_C$.
Then $\tor^p_{\mathcal{O}_C}(F_C,\ptfield(x_C)) \neq 0$ by the defining
property of $D_p$. We have (\cite{Mats86}, \S19, Lemma 1)
$$\prjdim_{\mathcal{O}_C} F_C = 
\sup \{ i \in \mathbb{Z} \; | \; \tor^i_{\mathcal{O}_C}(F_C,\ptfield(x_C))\}$$
hence $\prjdim_{\mathcal{O}_C} F_C \geq p$. By 
the Auslander-Buchsbaum equality we have 
$$\depth_{\mathcal{O}_C} \mathcal{O}_C = \prjdim_{\mathcal{O}_C} F_C +
\depth_{\mathcal{O}_C} F_C$$ and thus 
$ \codim_X C = \dim \mathcal{O}_C \geq \depth_{\mathcal{O}_C} \mathcal{O}_C
\geq p $ as required.
\end{proof}

The main idea behind the proof of the following proposition we owe 
to Bondal and Orlov in \cite{BonOr95}, Proposition 1.5.

\begin{proposition} \label{prps-spectral-argument}
Let $X$ be a locally noetherian scheme and $F \in D^b_{\text{coh}}(X)$ 
an object of finite $\tor$-dimension. Denote by 
$\mathcal{H}^i$ the $i$th cohomology sheaf of $F$. Then for any point 
$x \in X$ we have
\begin{align} \label{eqn-qmax-equal}
- \sup \{ i \in \mathbb{Z} \; | \; x \in \supp \mathcal{H}^i \}
= \inf \{ j \in \mathbb{Z} \; | \; \lder^j \iota_x^* F \neq 0 \}.
\end{align}

Let $C$ be an irreducible component of $\supp \mathcal{H}^{l}$
for some $l$ such that also $C \nsubseteq \supp \mathcal{H}^m$ 
for any $m < l$.
Then
\begin{align} \label{eqn-qmin-equal-1}
\codim_X C - \inf \{ i \in \mathbb{Z} \; | \; C \subseteq \supp
\mathcal{H}^i \}
= \sup \{ j \in \mathbb{Z} \; | \; \lder^j \iota^*_{x_C} F \neq 0 \}.
\end{align}
\end{proposition}
\begin{proof}
Fix a point $x \in X$. The main ingredient of the proof is 
the standard spectral sequence (eg. \cite{GelfandManin-MethodsOfHomologicalAlgebra}, Proposition III.7.10) associated 
to the filtration of $\lder\iota^*_x F$ by the rows of the 
Cartan-Eilenberg resolution of $F$:
\begin{align} \label{eqn-pullback-spectral-sequence}
E^{- p, q}_2 = \lder^{p} \iota_x^* (\mathcal{H}^q) \Rightarrow 
E^{q -p}_{\infty} = 
\lder^{p - q} \iota_x^*(F).
\end{align}

Denote by $h$ the highest non-zero row of $E^{\bullet \bullet}_2$. As all rows
above row $h$ and all columns to the right of column $0$ in 
$E^{\bullet\bullet}_2$ consist entirely of zeroes
\begin{center}
\includegraphics[scale=0.115]{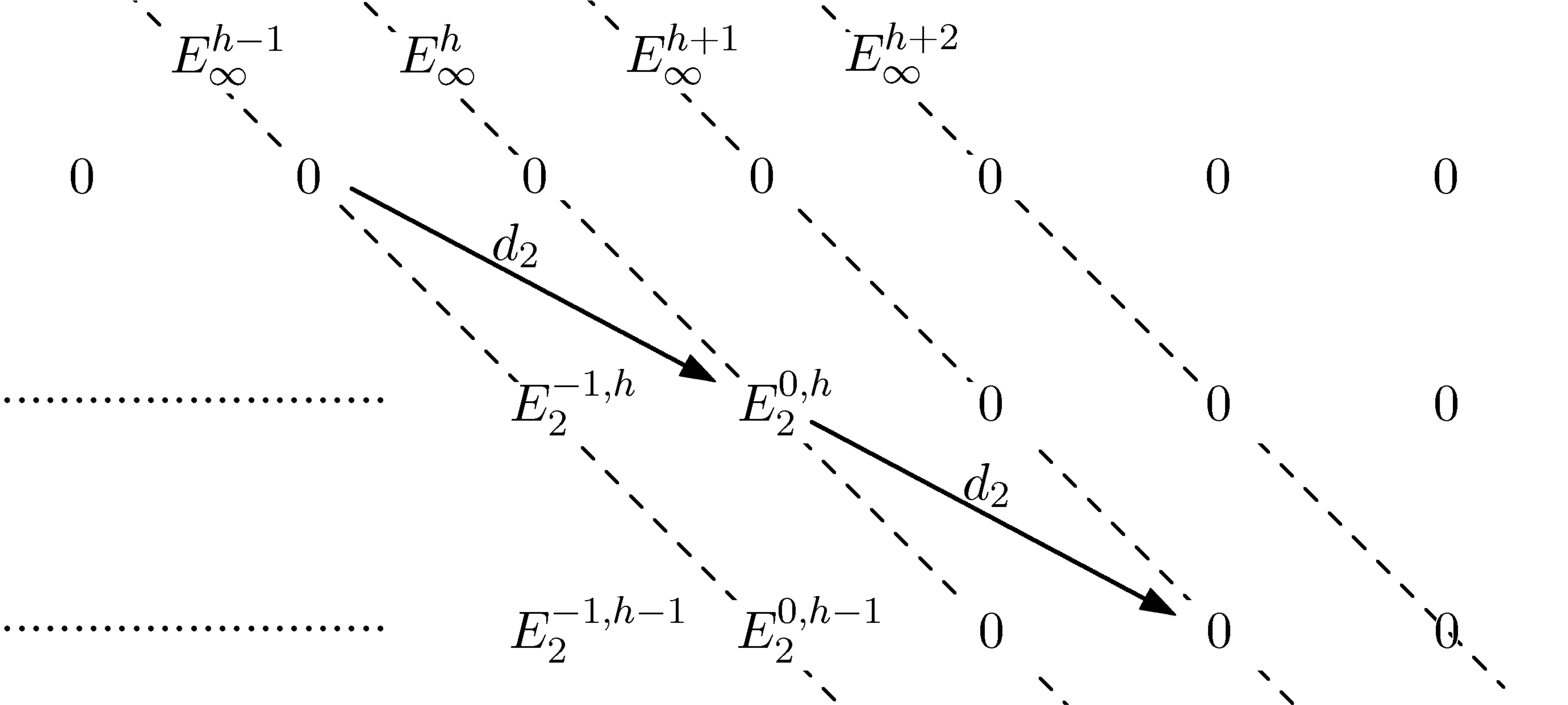}
\end{center}
we conclude by inspection of the complex that 
$0 = E^{n}_{\infty}$ for all $n > h$ and 
$\mathcal{H}^{h} |_x = E^{0,h}_2 = E^{h}_{\infty} =
\lder^{-h}(\iota^*_x(F))$. This gives $\eqref{eqn-qmax-equal}$.

To obtain $\eqref{eqn-qmin-equal-1}$ set $x$ to be the generic point of
$C$ and define $E^{\bullet \bullet}_{\bullet}$ as above.  
For any $m < l$ we have $C \nsubseteq \supp \mathcal{H}^m$ and 
hence $\lder \iota^*_x \mathcal{H}^m = 0$. So 
all the rows of $E^{\bullet \bullet}_2$ below $l$ consist of zeroes.
On the other hand, $C$ is an irreducible component of $\mathcal{H}^l$
and by Lemma $\ref{lemma-codim-of-cohdim-p-is-p}$ 
the set of points $y \in X$, such that there is a non-zero 
$\lder^i \iota^*_y (\mathcal{H}^l)$ with $i > d$, is closed and 
of codimension at least $d + 1$. Then this set can not contain $x$
for the closure of $x$ is $C$ whose codimension is $d$. Hence 
all columns to the left of column $-d$ in $E^{\bullet \bullet}_2$
consist entirely of zeroes. We conclude that 
$E^{n}_{\infty} = 0$ for all $n > l - d$ and 
$\lder^d \iota^*_x \mathcal{H}^l = E^{- d,l}_2 = E^{l - d}_\infty = 
\lder^{d-l} \iota^*_x F$. 
Thus, as $\lder^d \iota^*_x \mathcal{H}^l \neq 0$ 
by Lemma \ref{lemma-cohdim-at-least-support-codim}, we obtain
$\eqref{eqn-qmin-equal-1}$. 
\end{proof}

\begin{definition}
Let $\bf A\rm$ be an abelian category and $E^\bullet$ be a cochain complex
of objects of $\bf A\rm$. Define
its \it cohomological amplitude, \rm denoted by $\cxampl E^\bullet$,
to be the length of the minimal interval in $\mathbb{Z}$ containing
the set
\begin{align} \label{eqn-cohomological-dimension}
\{ i \in \mathbb{Z} \;|\; H^i(E^\bullet) \neq 0 \}.
\end{align} 
If no such interval exists we say that $\cxampl E = \infty$.
\end{definition}

Trivially $\cxampl E^\bullet$ is the minimal length of a bounded 
complex quasi-isomorphic to $E^\bullet$, if any exist, and 
infinity, if none do.  

\begin{definition}\label{defn-tor-amplitude}
Let $R$ be a ring or a sheaf of rings and $E^\bullet$ be 
a cochain complex of objects of $\modd$-$R$.  
Define its \it $\tor$-amplitude, \rm denoted by $\torampl_R E^\bullet$, 
to be the length of the minimal interval in $\mathbb{Z}$ 
containing the set 
\begin{align} \label{eqn-tor-amplitude}
\{ i \in \mathbb{Z} \;|\; \exists\; A \in \modd\text{-}R
\text{ such that }
\tor_R^i(E^\bullet, A) \neq 0 \}.
\end{align} 
If no such interval exists we say that $\torampl_R E = \infty$.
\end{definition}

Def. \ref{defn-tor-amplitude} can be seen to be equivalent to 
\cite{Kuznetsov-HomologicalProjectiveDuality}, Def. 2.20.

Let now $X$ be any scheme. It follows from 
\cite{Hartshorne-Residues-and-Duality}, Prop 4.2, that an object 
of $D^b(\modd\text{-}X)$ has finite $\tor$-amplitude if and only if it
is of finite $\tor$-dimension, i.e. quasi-isomorphic to 
a bounded complex of flat sheaves.
\begin{lemma} \label{lemma-pullbacks-to-tor-amplitude}
Let $X$ be a locally noetherian scheme and $E \in D^b_{\text{coh}}(X)$
an object of finite $\tor$-dimension. Denote by $l$ 
the length of the shortest complex of flat sheaves
quasi-isomorphic to $E$, and by $k$ the length of the smallest
interval in $\mathbb{Z}$ containing the set
\begin{align} \label{eqn-tor-amplitude-via-pullbacks}
\{ i \in \mathbb{Z} \;|\; \exists\; x \in X 
\text{ such that }
\lder^i \iota^*_x (E) \neq 0 \}.
\end{align} 
Then $l = \torampl_{\mathcal{O}_X} E = k$. 
\end{lemma}
\begin{proof}
Implications $ l \geq \torampl_{\mathcal{O}_X} E$ and
$\torampl_{\mathcal{O}_X} E \geq k$ are trivial. 
We claim that $k \geq l$. Let $n, k \in \mathbb{Z}$ be 
such that the interval $[-n-k,-n]$ contains the set 
$\eqref{eqn-tor-amplitude-via-pullbacks}$.  Then
$\eqref{eqn-qmax-equal}$ and $\eqref{eqn-qmin-equal-1}$
of Proposition $\ref{prps-spectral-argument}$ show that
$\mathcal{H}^i(E) = 0$ unless $i \in [n, n + k]$. Since resolutions by flat 
modules exist on $X$, there exists a complex $F^\bullet$ of flat 
sheaves quasi-isomorphic to $E$ and with $F_i = 0$ for all $i > n+k$. 
We claim that we can truncate $F^\bullet$ at degree $n$ and keep it flat, 
i.e.  that the sheaf $F^{n}/\img F^{n-1}$ is flat. But as
$\mathcal{H}^i(F^\bullet) = 0$ for $i < n$, the complex 
$$ \dots \rightarrow F^{n-2} \rightarrow F^{n-1} \rightarrow
F^{n} \rightarrow 0 \rightarrow  \dots $$
is a flat resolution of $F^{n}/\img F^{n-1}$. Hence $\lder^1 \iota^*_x 
(F^{n}/\img F^{n-1}) = \lder^{-n+1} \iota^*_x (E)$ and so 
vanishes for all $x \in X$ by assumption. Thus we obtain a
length $k$ complex of flat-sheaves quasi-isomorphic to $E$, i.e.
$k \geq l$.  
\end{proof}
Whenever $X$ is a quasi-projective scheme, or any other scheme where 
there exist resolutions by locally-free sheaves, replacing the word `flat' 
by the word `locally-free' throughout Lemma
\ref{lemma-pullbacks-to-tor-amplitude} and its proof shows that for any $E \in
D^b_{\text{coh}}(X)$ its $\tor$-amplitude is the length of the shortest complex 
of locally-free sheaves quasi-isomorphic to $E$. 
In other words, $\torampl_{\mathcal{O}_X} E$ is 
the \it homological dimension \rm 
of $E$ introduced in
\cite{BridgelandMacioca-FourierMukaiTransformsFor_K3_AndEllipticFibrations}.
The following can thus be compared to the inequality $ \text{hom.dim.} E \geq \codim C $
of
\cite{BridgelandMacioca-FourierMukaiTransformsFor_K3_AndEllipticFibrations}:

\begin{theorem}  \label{theorem-cohomological-amplitude-formula}
Let $X$ be a locally noetherian scheme and $E \in D^b_{\text{coh}}(X)$ 
an object of finite $\tor$-dimension. Then 
\begin{align} \label{eqn-global-cohomological-amplitude-formula}
\torampl_{\mathcal{O}_X} E \geq \codim \supp E + \cxampl E 
\end{align}
and for any irreducible component $C$ of $\supp E$ we have
\begin{align} \label{eqn-local-cohomological-amplitude-formula}
\torampl_{\mathcal{O}_C} E_C = \codim C + \cxampl_{\mathcal{O}_C} E_C.
\end{align}
\end{theorem}

\bf Remark: \rm To see that the inequality \eqref{eqn-global-cohomological-amplitude-formula} can be strict, consider 
$X = \mathbb{A}^1$ and $E = \mathcal{O}_X \oplus \mathcal{O}_x$ for
some closed point $x \in X$.

\begin{proof}
Denote by $\mathcal{H}^i$ the $i$th cohomology sheaf of $E$. Set
\begin{align*}
n = \inf_{x \in supp E} \{ i \in \mathbb{Z} | x \in \supp \mathcal{H}^i\} &&
m = \sup_{x \in supp E} \{ i \in \mathbb{Z} | x \in \supp \mathcal{H}^i\}
\\ 
l = \inf_{x \in supp E} \{ i \in \mathbb{Z} | \lder^i \iota_x^* E \neq 0 \} &&
h = \sup_{x \in supp E} \{ i \in \mathbb{Z} | \lder^i \iota_x^* E \neq 0 \} 
\end{align*}
and observe that $m - n = \cxampl E$ and $h - l =
\torampl_{\mathcal{O}_X} E$ (Lemma
\ref{lemma-pullbacks-to-tor-amplitude}). 

By $\eqref{eqn-qmax-equal}$ of Proposition \ref{prps-spectral-argument} 
we have 
\begin{align} \label{eqn-toramp-upperbound}
- m = l.
\end{align}
Let $D$ be any irreducible component of $\supp
\mathcal{H}^n$. We then have
\begin{align} \label{eqn-toramp-lowerbound}
\codim \supp E - n \leq \codim D - n = \sup \{ i \in \mathbb{Z} | \lder^i \iota_{x_D}^* E \neq 0 \}
\leq h
\end{align}
with the middle equality due to $\eqref{eqn-qmin-equal-1}$ of 
Proposition \ref{prps-spectral-argument} applied to $D$. Subtracting 
\eqref{eqn-toramp-upperbound} from \eqref{eqn-toramp-lowerbound} 
we obtain $(m - n) + \codim \supp E \leq (h - l)$. This shows
$\eqref{eqn-global-cohomological-amplitude-formula}$.

To obtain $\eqref{eqn-local-cohomological-amplitude-formula}$ we
observe that on $\spec \mathcal{O}_C$ the support of the localisation 
$E_C$ consists of a single point $x_C$. Therefore applying the above
argument to $X' = \spec \mathcal{O}_C$ and $E' = E_C$ we have 
$D = x_C = \supp E'$ which makes both the inequalities
in \eqref{eqn-toramp-lowerbound} into equalities.
\end{proof}

\section{Derived McKay correspondence} 
\label{section-main-proof}

Given a scheme $S$ denote by $D_{qc}(S)$ (resp. $D(S)$)  
the full subcategory of the derived category 
of $\mathcal{O}_S\text{-}\modd$ consisting of complexes 
with quasi-coherent (resp. bounded and coherent) cohomology. 
For $S$ a scheme of finite type over $\mathbb{C}$ and $H$ a finite
group acting on $S$ on the left by automorphisms an \it $H$-sheaf \rm
is a sheaf $\mathcal{E}$ of $\mathcal{O}_S$-modules equipped with a
lift of the $H$-action to $\mathcal{E}$. For the technical details see 
$\cite{BKR01}$, Section $4$. Denote by
$\mathcal{O}_S\text{-}\modd^H$ (resp. $\qcohcat^H S$, $\cohcat^H S$)
the abelian category of $H$-sheaves (resp. quasi-coherent, coherent
$H$-sheaves) on $S$ and by $D_{qc}^H(S)$ (resp. $D^H(S)$) 
the full subcategory of the derived category of 
$\mathcal{O}_S\text{-}\modd^H$ consisting of complexes with
quasi-coherent (resp. bounded and coherent) cohomology.

\subsection{Integral transforms}

\label{subsection-integral-transforms}

Let $N$ and $M$ be schemes of finite type over $\mathbb{C}$.
Denote by $\pi_N$ and $\pi_M$ the projections $N \times M \rightarrow
N$ and $N\times M \rightarrow M$. 
\begin{definition} \label{defn-integral-transform}
Let $E$ be an object 
of $D_{qc}(N \times M)$ of finite $\tor$-dimension. An 
\it integral transform $\Phi_E$ \rm is a functor 
$D_{qc}(N) \rightarrow D_{qc}(M)$ defined by
\begin{align} \label{eqn-integral-transform}
\Phi_E(-) = 
\rder \pi_{M *} (E \overset{\lder}{\otimes}
\pi^*_{N}(-)).
\end{align}
The object $E$ is called \it the kernel \rm of the transform.
If $\Phi_E$ is an equivalence of categories it is further
called \it a Fourier-Mukai transform\rm.
 
If a group $G$ acts on $N$ and $M$ then, for any $E \in D^G_{qc}(N
\times M)$ of finite $\tor$-dimension, \eqref{eqn-integral-transform}
defines an integral transform $D^G_{qc}(N) \rightarrow D^G_{qc}(M)$. If 
the group action on $N$ is trivial 
denote by $(- \otimes \rho_0)$ the functor $D_{qc}(N)~\rightarrow~
D^G_{qc}(N)$ which gives a sheaf the trivial $G$-equivariant 
structure. It is exact and has an exact left and right adjoint $(-)^G$,
the functor of taking the $G$-invariant part 
($\cite{BKR01}$, Section $4.2$). 
We also use the terms \it integral \rm and \it
Fourier-Mukai transform \rm for the functors $D_{qc}(N) \rightarrow
D^G_{qc}(M)$ of the form $\Phi_E(- \otimes \rho_0)$ where $\Phi_E$
is some integral transform $D^G_{qc}(N) \rightarrow D^G_{qc}(M)$.
\end{definition}

When $N$ and $M$ are smooth and proper varieties
it is well known that $\Phi_E$ has a left adjoint $\Phi_{E^\vee \otimes
\pi_M^*(\omega_M)[\dim M]}$ (\cite{BonOr95},
Lemma 1.2). The lemma below allows to generalise this 
to certain integral transforms between non-proper schemes. 
We use methods of Verdier-Deligne as per the exposition in \cite{DeligneCohomologieASupportPropre} 
to which we refer the reader for all the necessary definitions. 

\begin{lemma} \label{lemma-non-proper-adjoints}
Let $N$ and $M$ be separable schemes of finite type over $\mathbb{C}$
with $M$ smooth of dimension $n$. 
Let $E \in D(N \times M)$ be of finite homological dimension
with $\supp(E)$ proper over $N$.  
Then the functor 
$$\pi^*_N(-) \overset{\lder}{\otimes}
E:\quad D(N) \rightarrow D(N \times M)$$ has a left adjoint 
\begin{align} \label{eqn-nonproper-serre-duality}
\rder \pi_{N *}( - \overset{\lder}{\otimes} E^\vee \otimes
\pi^*_M(\omega_M))[n]:\quad D(N \times M) \rightarrow D(N).
\end{align}
\end{lemma}
\begin{proof}
First we compactify $M$: choose an open immersion 
$M \hookrightarrow \bar{M}$ with $\bar{M}$ smooth and proper 
\cite{Nagata-ImbeddingOfAnAbstractVarietyInACompleteVariety}. Then 
$\pi_N$ decomposes as an open immersion 
$\iota: N \times M \hookrightarrow N \times \bar{M}$ followed by 
the projection $\bar{\pi}_N: N \times \bar{M} \rightarrow N$. 
As $\bar{\pi}_N$ is smooth and proper Grothendieck-Serre
duality for smooth and proper morphisms  
(e.g. \cite{Hartshorne-Residues-and-Duality}, VII4.3) 
implies that $\bar{\pi}^*_{N}: D(N) \rightarrow D(N \times
\bar{M})$ has a left adjoint 
$$\rder \bar{\pi}_{N *}(-) \otimes \bar{\pi}^*_M \omega_{\bar{M}} [n]$$ 
where $\bar{\pi}_{M}:\; N \times \bar{M} \rightarrow
\bar{M}$ is the projection onto the second component. 

By the duality for open immersions
(\cite{DeligneCohomologieASupportPropre}, Prop. 4) 
 the left adjoint to the (exact) functor 
$\iota^*(-)$ exists as an (exact) functor $\iota_{!}$ from 
$\cohcat(N \times M)$ to the category $\proobjet$-$\cohcat(N \times
\bar{M})$. For the definition of $\proobjet$-$\cohcat(N \times
\bar{M})$ and the generalities on $\proobjet$-objects see
\cite{DeligneCohomologieASupportPropre}, $\text{n}^{\circ}$ 1. 
The functor $\iota_{!}$ may be calculated as follows: given 
$\mathcal{A} \in \cohcat(N \times M)$ take any $\bar{\mathcal{A}}
\in \cohcat(N \times \bar{M})$ which restricts to 
$\mathcal{A}$ on $N \times M$. Then 
\begin{align} \label{eqn-definition-of-prolongment-par-zero} 
\iota_!(\mathcal{A}) = \varinjlim \homm (\mathcal{I}^n
\bar{\mathcal{A}}, -)
\end{align}
where $\mathcal{I}$ is the ideal sheaf defining the complement
$N \times (\bar{M} \setminus M)$.

Finally, as $E$ is of finite homological dimension, 
the left adjoint of $(-)\overset{\lder}{\otimes}E$ is 
$(-) \overset{\lder}{\otimes} E^\vee$ where 
$E^\vee$ is 
$\rder \shhomm(E,\mathcal{O}_{N \times M})$. 

Therefore the left adjoint of $\pi^{*}_{N}(-) \overset{\lder}{\otimes}
E$ exists as the functor 
\begin{align} \label{eqn-equivariant-left-adjoint-proobjet}
\rder \bar{\pi}_{N *}(\iota_! (-\overset{\lder}{\otimes}
E^\vee) \otimes \bar{\pi}_{M}^*(\omega_{M}))[n]
\end{align}
from $\proobjet\text{-}D(N \times M)$ to $\proobjet\text{-}D(N)$. To 
finish the proof it suffices now to show that $\iota_!
(-\overset{\lder}{\otimes} E^\vee) =
\iota_*(-\overset{\lder}{\otimes} E^\vee)$. Then applying
the projection formula to $\iota_*(-\overset{\lder}{\otimes} E^\vee) \otimes
\bar{\pi}_{M}^*(\omega_{M})$ in 
$\eqref{eqn-equivariant-left-adjoint-proobjet}$ and observing that
$\iota \circ \bar{\pi}_{M} = \pi_M$ and $\iota \circ \bar{\pi}_{N} =
\pi_{N}$ yields $\eqref{eqn-nonproper-serre-duality}$.

We have $\id = \iota^* \iota_*$ on $\qcohcat(N \times
M)$ (\cite{Grothendieck-EGA-I}, \it Prop. 9.4.2\rm). 
It induces by the adjunction of \cite{DeligneCohomologieASupportPropre},
Prop. 4 natural transformations 
$\Upsilon:\; \iota_! \rightarrow \iota_*$  of functors 
$\cohcat(N \times M) \rightarrow \proobjet\text{-}\qcohcat(N \times
\bar{M})$ and $\Upsilon':\; \iota_!  (-\overset{\lder}{\otimes} E^\vee) \rightarrow 
\iota_*  (-\overset{\lder}{\otimes} E^\vee)$ of functors
$D(N\times M) \rightarrow \proobjet \text{-} D(N\times \bar{M})$.
By \cite{DeligneCohomologieASupportPropre}, Prop. 3 and 
the exactness of $\iota_!$ and $\iota_*$, to show $\Upsilon'$
to be an isomorphism of functors it suffices to show 
that $\Upsilon$ is an isomorphism on the cohomology sheaves of 
$-\overset{\lder}{\otimes} E^\vee$. The support
of these is proper over $N$ by the assumption on $E$.
For any $\mathcal{A} \in \cohcat(N \times M)$ 
we have 
\begin{align} \label{eqn-homs-shrek-pushdown}
\homm(\iota_!(\mathcal{A}), \iota_*(\mathcal{A})) = \varinjlim
\homm_{N \times \bar{M}}(\mathcal{I}^k \bar{\mathcal{A}}, 
\iota_*(\mathcal{A}))
\end{align}
using the notation of \eqref{eqn-definition-of-prolongment-par-zero}. 
From the construction of the adjunction in
\cite{DeligneCohomologieASupportPropre}, Prop. 4 it 
is immediate that $\Upsilon(\mathcal{A})$ is the unique element of
RHS of $\eqref{eqn-homs-shrek-pushdown}$ which restricts to $N \times
M$ as $\id \in \homm_{N \times M}(\mathcal{A},\mathcal{A})$. If
 $\supp(\mathcal{A})$ is proper over $N$, we can 
take $\bar{\mathcal{A}} = \iota_*\mathcal{A}$ in
\eqref{eqn-definition-of-prolongment-par-zero}. Moreover, $\mathcal{I}^k
\iota_*(\mathcal{A}) = \iota_*(\mathcal{A})$ for all $k$. Therefore
$\eqref{eqn-definition-of-prolongment-par-zero}$ yields
$\iota_!(\mathcal{A}) = \iota_*(\mathcal{A})$ and moreover
the RHS of \eqref{eqn-homs-shrek-pushdown} is just 
$\homm(\iota_* \mathcal{A}, \iota_* \mathcal{A})$. 
It is then clear that $\Upsilon(\mathcal{A})= \id$, as required.  
\end{proof}

\subsection{$G$-constellations and $\gnat$-families} 
\label{subsection-g-constellations}

\begin{definition}
Let $G$ be a finite subgroup of $\gl_n(\mathbb{C})$. 
A \it $G$-constellation \rm is a coherent $G$-sheaf $\mathcal{V}$ 
on $\mathbb{C}^n$ whose global sections $\Gamma(\mathcal{V})$ 
have the $G$-representation structure of the regular representation 
$\regrep$. 

Two $G$-constellations
$\mathcal{V}, \mathcal{W}$ are \em orthogonal 
in degree $k$ \rm if $G$-$\ext^k_{\mathbb{C}^n}(\mathcal{V},\mathcal{W}) = G$-$\ext^k_{\mathbb{C}^n}(\mathcal{W},\mathcal{V})= 0$. 
\end{definition}

Let now $Y$ be a scheme of finite type over $\mathbb{C}$. We endow
$Y$ with the trivial $G$-action, thus we can speak of $G$-sheaves
on $Y$ and on $Y \times \mathbb{C}^n$.  

\begin{definition} A \it $\gnat$-family on $Y$ \rm  
(short for \it$ G$-natural \rm or \it
geometrically natural\rm) is an object $\mathcal{F}$ 
of $\cohcat^G(Y \times \mathbb{C}^n)$, flat over $Y$, such that
for every closed $y \in Y$ the fiber $\mathcal{F}_{|y}$ is a
$G$-constellation supported on a single $G$-orbit. The 
\it Hilbert-Chow map \rm $\pi_\mathcal{F}$ of $\mathcal{F}$ is
the map $Y \rightarrow \mathbb{C}^n/G$ defined by 
$y \mapsto \supp_{\mathbb{C}^n} \mathcal{F}_{|y}$. A $\gnat$-family 
on a fixed morphism $Y \xrightarrow{\pi} \mathbb{C}^n/G$ is a
$\gnat$-family on $Y$ whose Hilbert-Chow map coincides with $\pi$.

Two subsets $C$ and $C'$ of $Y$ are \it orthogonal in degree $k$ 
in $\mathcal{F}$ \rm if for every $y \in C$ and $y' \in C'$ the fibers
$\mathcal{F}_{|{y}}$ and $\mathcal{F}_{|{y'}}$ are orthogonal in
degree $k$. The family $\mathcal{F}$ is \it orthogonal in degree $k$
\rm if $Y$ is orthogonal to $Y$ in degree $k$ in $\mathcal{F}$.
\end{definition}

\begin{proposition}\label{prps-hilbert-chow-map-is-a-morphism}
For any $\gnat$-family $\mathcal{F}$ its Hilbert-Chow 
map $\pi_\mathcal{F}$ is a morphism.
\end{proposition}
\begin{proof}
Denote by $\regring$ the ring $\mathbb{C}[x_1,\dots,x_n]$. 
For any $G$-constellation $\mathcal{V}$, 
the action of $\regring$ on $H^0(\mathcal{V})$ restricts
to the action of $\regring^G$ on $H^0(\mathcal{V})^G$. Clearly 
\begin{align} \label{eqn-pif-construction}
(\ann_\regring H^0(\mathcal{V}))^G \subseteq
\ann_{\regring^G} H^0(\mathcal{V})^G.
\end{align}
The LHS of $\eqref{eqn-pif-construction}$ is the image
of $\supp_{\mathbb{C}^n} \mathcal{V}$ in $\mathbb{C}^n/G$.
If this support is a single $G$-orbit, then 
$(\ann_\regring H^0(\mathcal{V}))^G$ is maximal in $\regring^G$
and $\eqref{eqn-pif-construction}$ is an equality.  
Therefore it suffices to construct a morphism $Y \rightarrow
\mathbb{C}^n/G$ which sends each $y \in Y$ to $\ann_{\regring^G}
H^0(\mathcal{F}_{|y})^G$. We construct it thus: the invariant part 
of $\pi_{Y*}(\mathcal{F})$ is a line bundle on $Y$, which
has a $\regring^G$-module 
structure induced from $\mathcal{F}$. This structure defines a homomorphism 
$\regring^G \rightarrow \mathcal{O}_Y$. The
corresponding morphism $Y \rightarrow \mathbb{C}^n/G$ 
is easily seen to send each $y \in Y$ to $\ann_{\regring^G}
H^0(\mathcal{F}_{|y})^G$.
\end{proof}

\begin{lemma}
\label{lemma-integral-transform-restricts-to-bounded-coherent}
If $\mathcal{F}$ is a $\gnat$-family on $Y$ and $\pi_\mathcal{F}: Y
\rightarrow \mathbb{C}^n/G$
is proper, then $\mathcal{F}$ is of finite homological dimension 
in $D^G(Y \times \mathbb{C}^n)$  and the integral transform $\Phi_\mathcal{F}: D_{qc}^G(Y) \rightarrow
D_{qc}^G(\mathbb{C}^n)$ restricts to   
$D^G(Y) \rightarrow D^G(\mathbb{C}^n)$.
\end{lemma}
\begin{proof}
Let $\iota$ be the open immersion $Y \times \mathbb{C}^n \rightarrow Y
\times \mathbb{P}^n$. As $\supp \mathcal{F}$ is proper over $Y$,
$\iota_* \mathcal{F}$ is coherent. Quite generally, given any coherent 
sheaf $\mathcal{A}$ on $Y \times \mathbb{P}^n$ flat over $Y$, 
consider the adjunction co-unit 
$\xi: \;\pi^*_Y \pi_{Y*} \mathcal{A} \rightarrow \mathcal{A}$. 
As $\pi_Y$ is proper and $\mathcal{A}$ is flat over $Y$,  
$\pi^*_Y \pi_{Y*} \mathcal{A}$ is lffr (locally free of finite rank). 
Twisting by some power of $\pi^*_{\mathbb{P}^n} \mathcal{O}(1)$
we can make $\xi$ surjective. But then $\ker \xi$ is again coherent 
and flat over $Y$. We set initially $\mathcal{A} = \iota_*
\mathcal{F}$ and repeat this construction until $\ker \xi$ becomes
lffr. This has to happen eventually as $\iota_* \mathcal{F}$ is flat
 over $Y$ and $\mathbb{P}^n$ is smooth. Thus we obtain an lffr
resolution of $\iota_*\mathcal{F}$ of finite length. Restricting it
to $Y \times \mathbb{C}^n$ demonstrates the first claim.
 
For the second claim: 
since $\pi_{Y}$ is flat, the pullback $\pi_{Y}^*(- \otimes \rho_0)$ is exact and
takes $D(Y)$ to $D^G(Y \times \mathbb{C}^n)$. 
Since $\mathcal{F}$ is of finite homological dimension,
$\mathcal{F} \overset{\lder}{\otimes} -$ takes $D^G(Y \times \mathbb{C}^n)$ to 
$D^G(Y \times \mathbb{C}^n)$. Moreover the image $\img(\mathcal{F}
\overset{\lder}{\otimes} -)$ lies in the full
subcategory of $D^G(Y \times \mathbb{C}^n)$ consisting of the objects
with support in $\supp \mathcal{F}$. Finally, 
$\pi_\mathcal{F}$ being proper implies that $\supp \mathcal{F}$ is
proper over $\mathbb{C}^n$, hence $\rder \pi_{\mathbb{C}^n *}$ takes 
$\img(\mathcal{F} \overset{\lder}{\otimes} -)$ to $D^G(\mathbb{C}^n)$ (\cite{Grothendieck-EGA-III-1}, \em Corollaire 3.2.4\rm). 
\end{proof}
 
The following demonstrates a certain relevance of \gnat-families: 

\begin{proposition}\label{prps-psps-necessary-criterion}
Let $G$ be a finite subgroup of $\gsl_n(\mathbb{C})$, $Y$
a variety and $E \in D^G(Y \times \mathbb{C}^n)$ an object such 
that $\Phi_{E}(- \otimes \rho_0)$ is an equivalence
$D(Y) \xrightarrow{\sim} D^G(\mathbb{C}^n)$ which sends
point sheaves on $Y$ to pure sheaves. Then $E$ is a \gnat-family over 
$Y$ and its Hilbert-Chow map $\pi_E$ is a crepant resolution of
$\mathbb{C}^n/G$. Moreover 
\begin{align} \label{eqn-simplicity-and-orthogonality}
G\text{-}\ext^i(E_{|y_1}, E_{|y_2}) = 
\begin{cases}
\mathbb{C} & \text{ if } y_1=y_2, i = 0 \\
0 & \text{ if } y_1 \neq y_2 
\end{cases}
\end{align}
and for any $y \in Y$ the (Kodaira-Spencer) map
$\ext^1(\mathcal{O}_y,\mathcal{O}_y) 
\rightarrow G\text{-}\ext^1(E_{|y},E_{|y})$ is an isomorphism.
\end{proposition}

\begin{proof}
By \cite{Huybrechts-FourierMukaiTransformsInAlgebraicGeometry}, 
Example 5.1(vi), $E_{|y} = \Phi_{E}(\mathcal{O}_y \otimes \rho_0)$,
whence the assertion $\eqref{eqn-simplicity-and-orthogonality}$
and the Kodaira-Spencer maps being isomorphisms. By \cite{Bridg97},
Lemma 4.3, it follows that $E$ is a pure sheaf flat over $Y$. 
Then by Lemma \ref{lemma-non-proper-adjoints} the inverse of 
$\Phi_{E}(- \otimes \rho_0)$ is $\Phi_{E^\vee[n]}(-)^G$. It maps $\mathcal{O}_{\mathbb{C}^n}$
to $(\pi_{Y *} E^\vee[n])^G$, so the cohomology sheaves of 
$(\pi_{Y *} E^\vee[n])^G$ are coherent $\mathcal{O}_Y$-modules.
Since $\pi_{Y *}$ is affine, the support of $E^\vee[n]$ is finite over $Y$. 
As $\supp(E^\vee[n]) = \supp
E$, we conclude that for each $y \in Y$ 
the support of $E_{|y}$ is a finite union of $G$-orbits. 
The simplicity of $E_{|y}$ further implies that it has 
to be a single $G$-orbit. To show that $\Gamma(E_{|y})$
has $G$-representation structure of $\regrep$ it suffices, 
by flatness of $E$, to show it for any single $y \in Y$. 
As the set $\{E_{|y}\}_{y \in Y}$ is an image of a spanning class
of $D(Y)$ under $\Phi(- \otimes \rho_0)$, it is a spanning class 
for $D^G(\mathbb{C}^n)$. Hence for every $G$-orbit $Z$ in 
$\mathbb{C}^n$ there exists $y \in Y$ such that $E_{|y}$ is supported 
at $Z$. Choose $Z$ to be any free orbit. The only simple 
$G$-sheaf supported on a free orbit is its structure sheaf, therefore 
$\Gamma(E_{|y}) \simeq \regrep$. We conclude that $E$ is a
gnat-family and that $\pi_E$ is surjective and an isomorphism outside 
the singularities of $\mathbb{C}^n/G$. By 
\cite{Roberts-MultiplicitiesAndChernClassesInLocalAlgebra}, Theorem 6.2.2
and \cite{BKR01}, Lemma 3.1, $Y$ is smooth and $\pi_E$ is crepant. 
It remains to show that $\pi_E$ is proper, which
is equivalent to $\supp_{Y \times \mathbb{C}^n} E$ being proper over $\mathbb{C}^n$ and that follows, e.g., from $\pi_{\mathbb{C}^n *} E$ 
having to be coherent, as it is a cohomology sheaf of the complex
$\Phi_E(\mathcal{O}_Y \otimes \rho_0)$.
\end{proof}

\subsection{Main results}

We now give the proof of Theorem \ref{theorem-main-theorem}. 
Its general framework follows those of \cite{BonOr95}, 
Theorem 1.1 and of \cite{BKR01}, Theorem 1.1. We note two
principal differences: \cite{BonOr95} works with smooth varieties, 
while we assume $Y$ to be a not necessarily smooth scheme (whence
the content of Section \ref{section-cohomological-width}); 
\cite{BKR01} adopts a two-step strategy to establish 
the left adjoint of $\Phi_\mathcal{F}(- \otimes \rho_0)$, whereas
our Lemma \ref{lemma-non-proper-adjoints} achieves this directly.

\begin{proof}[Proof of Theorem \ref{theorem-main-theorem}]

We divide the proof into five steps:

\it Step 1: We claim that $\Phi_\mathcal{F}(- \otimes \rho_0)$ has 
a left adjoint $(\Psi_{\mathcal{F}})^G$, where $\Psi_{\mathcal{F}}$ is
a certain integral transform $D^G(\mathbb{C}^n) \rightarrow D^G(Y)$. \rm

Recall that $\Phi_\mathcal{F} = \rder \pi_{\mathbb{C}^n *} (\mathcal{F}
\overset{\lder}{\otimes} \pi^*_Y (-))$. The issue here is the 
left adjoint of $\pi^*_{Y}(-)$ as $\pi_Y$, though smooth, is 
manifestly non-proper. But the support of $\mathcal{F}$ is proper,
so by Lemma \ref{lemma-non-proper-adjoints} the functor
$\rder \pi_{Y *} ( - \overset{\lder}{\otimes} \mathcal{F}^\vee[n])$ is the left adjoint
to $\pi^*_Y(-) \overset{\lder}{\otimes} \mathcal{F}$. The claim now
follows, for $\pi^*_{\mathbb{C}^n}$ is the left adjoint to 
$\rder \pi_{\mathbb{C}^n *}$ and $(-)^G$ is the left (and right) adjoint 
of $- \otimes \rho_0$.

\it Step 2: We claim that the composition $(\Psi_\mathcal{F})^G \circ
\Phi_\mathcal{F}(- \otimes \rho_0)$ is an integral transform $\Phi_Q$
for some $Q \in D(Y \times Y)$ and that for any closed point $(y_1, y_2)$ 
in $Y \times Y$ and any $k \in \mathbb{Z}$ we have 
\begin{align} \label{eqn-pullback-of-kernel-gives-gexts}
\lder^k \iota_{y_1, y_2}^* Q = G\text{-}\ext^k(\mathcal{F}_{| y_1},
\mathcal{F}_{| y_2})^\vee.
\end{align} \rm

The first assertion is a standard result due to Mukai 
in \cite{Mukai-DualityBetweenDXandDX^WithItsApplicationToPicardSheaves}, 
Proposition 1.3. The second assertion follows from 
the formula $(5)$ of $\cite{BKR01}$, Sec. 6, Step 2
by the adjunction of $\lder \iota^*_{y_1,y_2}$ and 
$\iota_{y_1,y_2*}$.

\it Step 3: We claim that $Q$ is a pure sheaf and that its support lies 
within the diagonal $Y \overset{\Delta}{\longrightarrow} Y \times Y$. \rm

First note that since $Y \times Y$ is of finite type over
$\mathbb{C}$, it is certainly Jacobson (see
\cite{Grothendieck-EGA-IV-3}, \S 10.3) and so any closed set of $Y
\times Y$ is uniquely identified by its set of closed points. We 
implicitly use this property at several points of the argument
below. 

Recall the closed set $N_k$ of $\eqref{eqn-codimension-of-nonorthogonality}$. 
As the support of any $G$-constellation is proper and as
$\omega_{\mathbb{C}^n} = \mathcal{O}_{\mathbb{C}^n} \otimes \rho_0$
as a $G$-sheaf since $G \subseteq \gsl_n(\mathbb{C})$, 
Serre duality applies to yield 
$$ G\text{-}\ext_{\mathbb{C}^n}^k(\mathcal{F}_{| y_1}, \mathcal{F}_{| y_2}) = 
   G\text{-}\ext_{\mathbb{C}^n}^{n-k}(\mathcal{F}_{| y_2},
\mathcal{F}_{| y_1})^\vee. $$
It follows that $\codim N_k = \codim N_{n - k}$ for all $k$.

Let $C$ be an irreducible component of $\supp Q$. Denote by $y_C$ its
generic point, by $\mathcal{O}_C$ the local ring of $y_C$
and by $Q_C$ the localisation of $Q$ to $\mathcal{O}_C$.  
For any $k$ denote by $M_k$ the set 
$\{y \in Y \times Y \;|\; \lder^k \iota_y^*Q \neq 0\}$
and let $l$ and $m$ be the infimum and the supremum of the set
$\{k \in \mathbb{Z} \;|\; y_C \in M_k \}$, hence 
$\torampl_{\mathcal{O}_C} Q_C = m - l $ (Lemma
\ref{lemma-pullbacks-to-tor-amplitude}). 
By $\eqref{eqn-pullback-of-kernel-gives-gexts}$ 
the closure of $M_l \setminus \Delta$ is $N_l$, so
$y_C \in M_l$ implies $y_C \in \Delta$ or $y_C \in N_l$.
Similarly for $N_m$. Thus either $y_C \in \Delta$ or $y_C \in N_l \cap
N_m$. The latter would imply that 
\begin{align*}
&\codim C \geq \codim N_l \geq n - 2l + 1 \\
&\codim C \geq \codim N_m = \codim N_{n-m} \geq 2m - n + 1 
\end{align*}
and therefore that $\;\codim C \geq m - l + 1$. But then 
$\;\codim C $ would be strictly greater than $\torampl_{\mathcal{O}_C}
Q_C$, which contradicts 
Theorem \ref{theorem-cohomological-amplitude-formula}. Thus $y_C$ 
lies within $\Delta$ and, since $Y$ is separated, so does all of $C$.  

We have now shown that $\supp Q \subseteq \Delta$, so $\codim \supp Q \geq
n$. But as $\mathbb{C}^n$ is smooth and $n$-dimensional, 
$\eqref{eqn-pullback-of-kernel-gives-gexts}$ implies 
\begin{align} \label{eqn-pullbacks-vanish-outside-0-to-n}
& \lder^k \iota^*_y Q = 0 & &\forall y \in Y,\; k \notin 0, \dots, n 
\end{align}
so $\torampl Q \leq n$. 
By Theorem $\ref{theorem-cohomological-amplitude-formula}$ 
$\torampl Q = n$ and
$\cxampl Q = 0$. Together with 
$\eqref{eqn-pullbacks-vanish-outside-0-to-n}$ this 
implies that $Q$ is a pure sheaf. 

\it Step 4: We claim that $Q$ is the structure sheaf
$\mathcal{O}_\Delta$ of the diagonal $\Delta$ and
therefore $\Phi_{\mathcal{F}}(- \otimes \rho_0)$ is 
fully faithful. \rm 

The adjunction co-unit $\Phi_Q \rightarrow \id_{D(Y)}$ induces
a surjective $\mathcal{O}_{Y\times Y}$-module morphism $Q \xrightarrow{\epsilon}
\mathcal{O}_{\Delta}$. Let $K$ be its kernel, we then
have a short exact sequence
\begin{align} \label{eqn-Q-O-Delta-ses}
0 \rightarrow K \rightarrow Q \xrightarrow{\epsilon} \mathcal{O}_{\Delta}
\rightarrow 0.
\end{align}
Choosing some closed point $(y,y) \in \Delta$ and applying
functor $\lder \iota_{y,y}^* (-)$ to \eqref{eqn-Q-O-Delta-ses}
we obtain a long exact sequence of $\mathbb{C}$-modules
\begin{align*}
\dots \rightarrow G\text{-}\ext^{1}_{\mathbb{C}^n}(\mathcal{F}_{|y},
\mathcal{F}_{|y})^* \xrightarrow{\alpha_y} \Omega^1_{Y,y} \rightarrow 
K_{y,y} \rightarrow G\text{-}\eend_{\mathbb{C}^n}(\mathcal{F}_{|y})^* \xrightarrow{\epsilon_y} \mathbb{C} 
\rightarrow 0 \rightarrow \dots .
\end{align*}
The map $\epsilon_y$ is surjective due to any $G$-constellation having
automorphisms consisting of scalar multiplication. 
It is an isomorphism whenever $\mathcal{F}_{|y}$ is simple, 
i.e. when the scalar multiplication automorphisms are all we get. 
The map $\alpha_y$ is the dual of the Kodaira-Spencer map of
$\mathcal{F}$ at $y \in Y$, which takes a tangent vector at $y$ 
to the infinitesimal deformation in that direction in the family 
$\mathcal{F}$. Hence 
for any $y \in Y$, such that $\mathcal{F}_{|y}$ is simple and 
such that the Kodaira-Spencer map of $\mathcal{F}$ is injective at $y$, 
the long exact sequence above shows that $K|_{y,y} = 0$.  

Having proved that $\supp Q \subseteq \Delta$ we have proved
by \eqref{eqn-pullback-of-kernel-gives-gexts} that any 
two $G$-constellations in $\mathcal{F}$ are orthogonal. 
Denoting by $q$ the quotient map $\mathbb{C}^n
\rightarrow \mathbb{C}^n/G$ we claim that for any closed point 
$x \in \mathbb{C}^n/G$, such that $q^{-1}(x)$ is a free orbit of 
$G$, the fiber $\pi_{\mathcal{F}}^{-1}(x)$ consists of at most 
a single point. This is because, by definition of
$\pi_{\mathcal{F}}$, all the $G$-constellations parametrised by 
$\pi_{\mathcal{F}}^{-1}(x)$ are supported on $q^{-1}(x)$ - 
and any two $G$-constellations supported 
at the same free orbit are easily seen to be isomorphic. 
Thus $\pi_\mathcal{F}$ is an isomorphism on the smooth locus $X_0$ 
of $\mathbb{C}^n/G$. 
By $\cite{Logvinenko-Natural-G-Constellation-Families}$, Proposition 1.5
the family $\mathcal{F}$ on $X_0$ (identified with an open subset of $Y$ 
via $\pi_\mathcal{F}$) is locally isomorphic to the canonical 
$G$-cluster family $q_* \mathcal{O}_{\mathbb{C}^n}|_{X_0}$. 
As any $G$-cluster is simple and as the Kodaira-Spencer map 
of $q_* \mathcal{O}_{\mathbb{C}^n}|_{X_0}$ is trivially injective 
$K|_{y,y} = 0$ for any $y \in X_0$. Therefore 
$\codim_{Y \times Y} \supp K \geq n+1$, as $X_0$ is open in $\Delta$. 
 
On the other hand, since $\torampl Q =
\torampl \mathcal{O}_\Delta = n$,
the short exact sequence \eqref{eqn-Q-O-Delta-ses} implies 
that $\torampl K \leq n$. As that is smaller than the
codimension of its support, $K = 0$ by
Theorem $\ref{theorem-cohomological-amplitude-formula}$. 
Thus $Q \simeq \mathcal{O}_\Delta$, the adjunction 
co-unit is an isomorphism and $\Phi_{\mathcal{F}}(- \otimes
\rho_0)$ is fully faithful. 

\it Step 5: We claim that $\Phi_{\mathcal{F}}(- \otimes \rho_0)$ is 
an equivalence of categories. \rm

As $D(Y)$ is fully faithfully embedded in $D^G(\mathbb{C}^n)$ 
the trivial Serre functor of the latter induces a trivial Serre functor
on the former. Therefore the left adjoint to 
$\Phi_{\mathcal{F}}(- \otimes \rho_0)$ is also its right adjoint. 
Then $\Phi_{\mathcal{F}}(- \otimes \rho_0)$ is an equivalence
of categories by $\cite{Bridg97}$, Theorem 3.3.

\end{proof}

\begin{proof}[Proof of Corollary \ref{cor-3d-application-of-main-theorem}]

It suffices to demonstrate that $\mathcal{F}$ satisfies the condition 
of Theorem $\ref{theorem-main-theorem}$. Thus we have to show that 
$\codim N_0 \geq 4$ and $\codim N_1 \geq 2$. 
But, as seen in the proof of Theorem $\ref{theorem-main-theorem}$, 
$N_k$ lies within the fibre product $Y \times_{\mathbb{C}^3/G} Y$ for all 
$k$. As $\pi_\mathcal{F}$ is birational its fibres are at most divisors 
and so the codimension of $Y \times_{\mathbb{C}^3/G} Y$ is at least $2$. 

It remains to show that $N_0 \geq 4$. The assumptions of the Corollary
ensure that $N_0$ is contained
in the union of all sets of form $(E_i \cap E_j) \times 
(E_k \cap E_l)$ or $E_i \times (E_i \cap E_j \cap E_k)$, 
and the codimension of each of these sets is $4$.  
\end{proof}

\section{Orthogonality in degree zero} \label{section-degree-zero-orthogonality}

Throughout this section we denote by $G$ a finite abelian subgroup 
of $\gsl_n(\mathbb{C})$, by $Y$ a smooth scheme of finite type over 
$\mathbb{C}$ and by $\mathcal{F}$ a $\gnat$-family on $Y$.
We assume that the Hilbert-Chow morphism $\pi_\mathcal{F}$ associated to
$\mathcal{F}$ is birational and proper. 
The main purpose of this section is to show how, given any pair of 
closed points of $Y$, one checks whether the corresponding pair of
$G$-constellations are orthogonal in degree $0$.

We denote by $\givrep$ the representation
of $G$ given by its inclusion into $\gsl_n(\mathbb{C})$. The (left) action 
of $G$ on $\givrep$ induces a right action of $G$ on $\givrep^\vee$ which 
we make into a left action by setting:
\begin{align} \label{eq-gaction}
  g\cdot f(v) = f(g^{-1} \cdot v)\quad\quad\text{ for all } 
v \in \givrep, \; f \in \givrep^\vee, \; g \in G.
\end{align} 
We denote by $x_1, \dots, x_n$ the common eigenvectors of the action
of $G$ on $\givrep^\vee$.  We denote by
$\regring$ the symmetric algebra $S(\givrep^\vee)$ with the induced
left action of $G$. Then $\regring =
\mathbb{C}[x_1,\dots,x_n]$ and as an affine $G$-scheme $\mathbb{C}^n$ 
is $\spec \regring$. We denote by $G^\vee$ the character group 
$\homm(G, \mathbb{C}^*)$ of $G$. A rational function $f \in K(\mathbb{C}^n)$ 
is said to be $G$-homogeneous of weight $\chi \in G^\vee$ if we have
$f(g.v) = \chi(g)\; f(v)$ for all $v \in \mathbb{C}^n$ where
$f$ is defined. We denote by $\rho(f)$ the weight $\chi$ of such $f$. 
It follows from \eqref{eq-gaction} that $G$ acts on $f$ by $\rho(f)^{-1}$.  

From here on we employ freely the terminology and the 
results of $\cite{Logvinenko-Natural-G-Constellation-Families}$.

\subsection{The McKay quiver of $G$}
\label{section-the-mckay-quiver-of-G}

By a \it quiver \rm we mean a vertex set $Q_0$, an arrow set $Q_1$ 
and a pair of maps $h \colon Q_1 \rightarrow Q_0$ and 
$t \colon Q_1 \rightarrow Q_0$ giving the head $hq \in Q_0$ and the
tail $tq \in Q_0$ of each arrow $q \in Q_1$. By a \it representation
of a quiver \rm we mean a graded vector space 
$\bigoplus_{i \in Q_0} V_i$ and a collection of linear maps
$\{\alpha_q\colon V_{tq} \rightarrow V_{hq}\}_{q \in Q_1}$. 

\begin{definition}
The \it McKay quiver of $G$ \rm is the quiver whose vertex set $Q_0$ 
are the irreducible representations $\rho$ of $G$ and whose arrow set $Q_1$ 
has $\dim \homm_G (\rho_i, \rho_j \otimes \givrep)$ arrows
going from the vertex $\rho_i$ to the vertex $\rho_j$.
\end{definition}

We have $\givrep^\vee = \bigoplus \mathbb{C} x_i$, as
$G$-representations. Denote by $U_\chi$ the $1$-dimensional representation 
on which $G$ acts by $\chi \in G^\vee$. By Schur's lemma 
\begin{align*}
G\text{-}\homm(U_{\chi_i} \otimes \givrep^{\vee}, U_{\chi_j}) =  
\begin{cases}
\mathbb{C} \quad \text{ if } \chi_j = \chi_i \rho(x_k)^{-1} \quad 
k \in \{1,\dots,n\} \\
0 \quad \text{ otherwise }
\end{cases}.
\end{align*}
Thus each vertex $\chi$ of the McKay quiver of $G$ has $n$ arrows
emerging from it and going to vertices $\chi \rho(x_k)^{-1}$ for $k =
1,\dots,n$. We denote the arrow from $\chi$ to $\chi \rho(x_k)^{-1}$ by 
$(\chi, x_k)$. Let now $A$ be a $G$-constellation viewed as an 
$\twalg$-module 
($\cite{Logvinenko-Natural-G-Constellation-Families}$, Section 1.1)
and let $\oplus A_\chi$ be its decomposition into irreducible
representations of $G$. Then the $\twalg$-module structure on $A$
defines a representation of the McKay quiver into the graded 
vector space $\oplus A_\chi$, where the map $\alpha_{\chi,
x_k}$ is just the multiplication by $x_k$, i.e.  
\begin{align} \label{eqn-gcon-rep-definition}
  \alpha_{\chi, x_k}:\; A_\chi \rightarrow A_{\chi\rho(x_k)^{-1}}, \; 
   v \mapsto x_k \cdot v.
\end{align}

\subsection{Degree $0$ orthogonality of $G$-constellations}
\label{subsection-degree-0-orthogonality-of-G-constellations}

Let $A$ and $A'$ be two $G$-constellations and $\phi$ be an
$\twalg$-module morphism $A \rightarrow A'$. Let $
\bigoplus_{G^\vee} A_\chi$ and $\bigoplus_{G^\vee} A'_\chi$ be
decompositions of $A$ and $A'$ into one-dimensional representations
of $G$. By $G$-equivariance $\phi$ decomposes into linear maps 
$\phi_\chi:\; A_\chi \rightarrow A'_\chi$. 

Let $\{ \alpha_q \}$ and $\{ \alpha'_q\}$ be the corresponding 
representations of the McKay quiver into graded vector spaces 
$\oplus A_\chi$ and $\oplus A'_\chi$, as per 
$\eqref{eqn-gcon-rep-definition}$. Each $\alpha_q$ is a linear map 
between one-dimensional vector spaces $A_{tq}$ and $A_{hq}$ and so 
is either a zero-map or an isomorphism, similarly for the maps 
$\alpha'_q$. So for each arrow of the McKay quiver we distinguish the 
following four possibilities:  

\begin{definition} \label{defn-arrow-type}
Let $q$ be an arrow of McKay quiver of $G$. With the notation
above we say that with respect to an ordered pair $(A,A')$ of 
$G$-constellations the arrow $q$ is: 
\begin{enumerate}
\item a type $[1,1]$ arrow, if both $\alpha_q$ and $\alpha'_q$ are 
isomorphisms.
\item a type $[1,0]$ arrow, if $\alpha_q$ is an isomorphism and 
$\alpha'_q$ is a zero map. 
\item a type $[0,1]$ arrow, if $\alpha_q$ is a zero map and 
$\alpha'_q$ is an isomorphism. 
\item a type $[0,0]$ arrow, if both $\alpha_q$ and $\alpha'_q$ are 
zero maps.
\end{enumerate}
\end{definition}

\begin{proposition}
Let $q$ and $(A,A')$ be as in Definition $\ref{defn-arrow-type}$
and let $\phi$ be any $\twalg$-module morphism $A \rightarrow A'$.
Then:
\begin{enumerate}
\item If $q$ is a $[1,0]$ arrow, then $A_{hq} \subseteq \ker \phi$.
\item If $q$ is a $[0,1]$ arrow, then $A_{tq} \subseteq \ker \phi$.
\item If $q$ is a $[1,1]$ arrow, $A_{tq}$ and $A_{hq}$ either
both lie in $\ker \phi$ or both don't. 
\end{enumerate}
\end{proposition}
\begin{proof}
Write $q = (\chi, i)$ where $\chi \in G^\vee$ and $i \in \{1,\dots,n\}$.
Recall that $\alpha_q$ is the map $A_{tq} \rightarrow A_{hq}$ 
corresponding to the action of $x_i$ on $A_{tq}$. Then $R$-equivariance 
of the morphism $\phi$ implies a commutative square
\begin{align*}
\xymatrix{
A_{hq} \ar[r]^{\phi_{hq}} & A'_{hq} \\
A_{tq} \ar[r]_{\phi_{tq}} \ar[u]^{\alpha_q} & A'_{tq} \ar[u]_{\alpha'_q}
}
\end{align*}
from which all three claims immediately follow.
\end{proof}

\begin{corollary} \label{cor-orthogonality-criterion}
Let $(A,A')$ be an ordered pair of $G$-constellations. If every 
component of the McKay quiver path-connected by $[1,1]$-arrows
has either a $[0,1]$-arrow emerging from it or a $[1,0]$-arrow 
entering it, then 
$$ \homm_{\twalg}(A,A') = 0.$$
If, also, every component has either a $[0,1]$-arrow entering it 
or a $[1,0]$-arrow emerging from it, then we further have
$$ \homm_{\twalg}(A',A) = 0 $$
and therefore $A$ and $A'$ are orthogonal in degree $0$. 
\end{corollary}

\subsection{Divisors of zeroes}

The Hilbert-Chow morphism $\pi_\mathcal{F}: Y \rightarrow
\mathbb{C}^n/G$ is birational, thus it defines a notion of
$G$-Cartier and $G$-Weil divisors on $Y$ 
($\cite{Logvinenko-Natural-G-Constellation-Families}$), Section 2). 
The family $\mathcal{F}$, in a sense of a sheaf of $\mathcal{O}_Y
\otimes (\twalg)$-modules on $Y$, can be written as 
$\bigoplus_{\chi \in G^\vee} \mathcal{L}(-D_\chi)$, where 
$D_\chi$ are $G$-Weil divisors. For any other such 
expression $\bigoplus \mathcal{L}(-D'_\chi)$ of $\mathcal{F}$ 
there exist $f \in K(Y)$ such that $D'_\chi = D_\chi + (f)$ for all 
$\chi \in G^\vee$ ($\cite{Logvinenko-Natural-G-Constellation-Families}$, Section 3.1).

\begin{definition} Let $q = (\chi, x_k)$ be an arrow in the
McKay quiver of $G$. We define the \it divisor of zeroes \rm
$B_{q}$ of $q$ in $\mathcal{F}$ to be the Weil divisor
\begin{align} \label{eqn-divisor-of-zeroes}
D_{\chi^{-1}} + (x_i) - D_{\chi^{-1}\rho(x_i)}.
\end{align}
\end{definition}

Note that $B_{q}$ is always an ordinary, integral Weil divisor on $Y$.  

\begin{proposition} \label{prps-divisors-of-zeroes}
Let $(\chi, x_k)$ be an arrow in the McKay quiver of $G$ and 
$B_{\chi, x_k}$ be its divisor of zeroes in $\mathcal{F}$. 
Let $y$ be a closed point of $Y$ and $A$ be the 
$G$-constellation $\mathcal{F}_{|y}$. Then in the corresponding 
representation $\{\alpha_q\}_{q \in Q_1}$ of 
the McKay quiver the map $\alpha_{\chi,x_k}$ is a zero map if and only 
if $y \in B_{\chi, x_k}$.
\end{proposition}
\begin{proof} The map $\alpha_{\chi,x_k}:\; A_\chi \rightarrow A_{\chi\rho(x_k)^{-1}}$ 
is the action of $x_k$ on $A_\chi$. This map 
is the restriction to the point $y$ of the global section 
$\beta$ of the $\mathcal{O}_Y$-module 
\begin{align}\label{eqn-action-map-section}
\shhomm_{G,\mathcal{O}_Y}(\mathcal{O}_Y x_k \otimes \mathcal{F}_\chi,
\mathcal{F}_{\chi\rho^{-1}(x_k)})
\end{align}
defined by $x_k \otimes s \mapsto x_k \cdot s$ for any section $s$ 
of the $\chi$-eigensheaf $\mathcal{F_\chi}$. 

As $G$ acts on a monomial of weight $\chi$ by $\chi^{-1}$ 
the $\chi$-eigensheaf of $\mathcal{F}$ 
is $\mathcal{L}(-D_{\chi^{-1}})$. Hence
$\eqref{eqn-action-map-section}$ is canonically isomorphic to 
the following sub-$\mathcal{O}_Y$-module of $K(\mathbb{C}^n)$:
\begin{align}\label{eqn-action-map-oy-module}
\mathcal{L}(D_{\chi^{-1}} + (x_k) - D_{\chi^{-1}\rho(x_k)})
\end{align}
and the isomorphism maps $\beta$ to the global section 
$1 \in K(\mathbb{C}^n)$ of $\eqref{eqn-action-map-oy-module}$.
Which vanishes precisely on the Weil divisor 
$B_{\chi, x_k} = D_{\chi^{-1}} + (x_k) - D_{\chi^{-1}\rho(x_k)}$.
\end{proof}

Proposition $\ref{prps-divisors-of-zeroes}$ together with Corollary 
$\ref{cor-orthogonality-criterion}$ show that the data of the 
divisors of zeroes of $\mathcal{F}$ is all that is necessary to 
determine whether any given pair of closed points of $Y$ are orthogonal 
in degree $0$ in $\mathcal{F}$.

\subsection{Direct transforms}
\label{section-direct-transforms-and-stable-families}

Let $Y'$ and $Y''$ be two crepant resolutions of $\mathbb{C}^n/G$ 
isomorphic outside of a closed set of codimension $\geq 2$. 
E.g. the case $n=3$ where all crepant 
resolutions are related by a chain of flops (\cite{Kollar-Flops}). We 
fix a birational isomorphism and use it to identify $Y'$ and $Y''$
along the isomorphism locus $U$. Since the complement of $U$ is 
of codimension $\geq 2$ in $Y'$ (resp. $Y''$) any line bundle or divisor 
on $U$ extends uniquely to a line bundle or a divisor on $Y'$ 
(resp. $Y''$). The same is true of 
a family of $G$-constellations as for $G$ abelian any such family 
is a direct sum of line bundles. For any family $\mathcal{V}'$ of 
$G$-constellations on $Y'$ we define its \it direct transform \rm
$\mathcal{V}''$ to $Y''$ to 
be the unique extension to $Y''$ of the restriction of $\mathcal{V}'$ to 
$U$. Observe that if $\mathcal{V}'$ is of form 
$\bigoplus_\chi \mathcal{L}(-D'_\chi)$ for some $G$-Weil divisors
$D'_\chi$ on $Y'$ then $\mathcal{V}''$ is the family 
$\bigoplus \mathcal{L}(-D''_\chi)$ where each $D''_\chi$ is the direct 
transform of $D'_\chi$. 

If $\mathcal{F}$ can be shown to be a direct transform 
of some everywhere orthogonal in degree $0$ family $\mathcal{F}'$ 
on some $Y'$, it greatly reduces the number of calculations necessary 
to determine the degree $0$ non-orthogonality locus of $\mathcal{F}$. 
Let $U$ be as above. 
As $\mathcal{F}$ is the direct transform of $F'$, the restriction of
$\mathcal{F}$ to $U \subset Y$ is isomorphic to the restriction of
$F'$ to $U \subset Y'$. So the calculations only have
to be carried out for points in $Y \times Y \setminus U \times U$. 

\subsection{Theta stability and $\gnat$-families}
\label{subsection-theta-stability}

We recall basic facts about $\theta$-stability for $G$-constellations, 
cf. \cite{Craw-Ishii-02}, Section 2.1. Let 
$\mathbb{Z}(G) = \bigoplus_{\chi \in G^\vee} \mathbb{Z} \chi$ 
be the representation ring of $G$ and set 
$$ \Theta = \{ \theta \in \homm_\mathbb{Z}(\mathbb{Z}(G),\mathbb{Q}) \;|\; 
\theta(\regrep) = 0 \} $$
For any $\theta \in \Theta$, a $G$-constellation $A$ is 
$\theta$-stable (resp. 
$\theta$-semistable) if for every sub-$\twalg$-module $B$ of $A$ we have
$\theta(B) > 0$ (resp. $\theta(B) \geq 0$). We say that $\theta$ is generic 
if every $\theta$-semistable $G$-constellation is $\theta$-stable. This 
is equivalent to $\theta$ being non-zero on any proper subrepresentation 
of $\regrep$. 

Let $\pi$ be any proper birational morphism $Y \rightarrow
\mathbb{C}^n/G$. A $\gnat$-family $\mathcal{V}$ 
on $Y \xrightarrow{\pi}
\mathbb{C}^n/G$ is \it normalized \rm 
if $\mathcal{V}^G \simeq \mathcal{O}_Y$. Such $\mathcal{V}$ 
can be written uniquely as $\bigoplus_{\chi \in G^\vee} 
\mathcal{L}(-D_\chi)$ for some $G$-Weil divisors $D_\chi$ with 
$D_\trch_ = 0$ 
($\cite{Logvinenko-Natural-G-Constellation-Families}$, Cor. 3.5). 
Denote by $\mathfrak{E}$ the set of all prime Weil divisors on $Y$
whose image in $\mathbb{C}^n/G$ is either a point or a coordinate
hyperplane $x_i^{|G|} = 0$. As $G$ is
abelian, any branch divisor of $\mathbb{C}^n \rightarrow
\mathbb{C}^n/G$, if it exists, is one of the hyperplanes $x_i^{|G|} = 0$.
Hence, by $\cite{Logvinenko-Natural-G-Constellation-Families}$, Prop. 3.14
and 3.15, each $D_\chi$ is of form $\sum_{E \in \mathfrak{E}}
q_{\chi,E} E$. Denote by $U$ the open subset of $Y$ consisting of 
points lying on at most one divisor in $\mathfrak{E}$. 

\begin{definition}
Let $\theta$ be an element of $\Theta$. 
We define a map $$w_\theta:\quad \left\{
\text{normalized $\gnat$-families on $Y\xrightarrow{\pi}
\mathbb{C}^n/G$} \right\} \rightarrow \mathbb{Q}$$ by
\begin{align}
w_\theta(\mathcal{V}) = \sum_{E \in \mathfrak{E}} \sum_{\chi \in
G^\vee} \theta(\chi) q_{\chi,E}.
\end{align}
\end{definition}

The domain of definition of $w_\theta$ is finite 
($\cite{Logvinenko-Natural-G-Constellation-Families}$, Corollary
3.16), so for any $\theta \in \Theta$ there is at least one normalized 
$\gnat$-family maximizing $w_\theta$. 
\begin{proposition}\label{prps-stability-criterion}
Let $\mathcal{M}$ be any family which maximizes $w_\theta(\mathcal{M})$. 
Then for any point $y \in U$ the fiber of $\mathcal{M}$ at $y$ 
is a $\theta$-semistable $G$-constellation. If, moreover, $\theta$ is 
generic, then such family $\mathcal{M}$ is unique. 
\end{proposition}
\begin{proof}
Write $\mathcal{M}$ as $\bigoplus \mathcal{L}(-M_\chi)$.
Suppose that the fiber of $\mathcal{M}$ is not $\theta$-semistable at some 
$y \in U$. Denote this fiber by $A$, its decomposition into
irreducible representations by $\bigoplus_{\chi \in G^\vee} A_\chi$
and the corresponding representation of the McKay quiver by
$\{\alpha_q\}$. As $A$ isn't $\theta$-semistable
 there exists a non-empty proper subset $I$ 
of $G^\vee$ such that $A' = \bigoplus_{\chi \in I} A_\chi$ is a 
sub-$\twalg$-module of $A$ and $\theta(A') < 0$. Denote by $J$
the complement $G^\vee \setminus I$. Denote by $Q_{I \rightarrow J}$ 
the subset $\{ q \in Q_1 \;|\; tq \in I, hq \in J \}$ of the arrow set 
$Q_1$ of the McKay quiver and 
similarly for $Q_{J \rightarrow I}$, $Q_{I \rightarrow I}$, $Q_{J
\rightarrow J}$. Then $A'$ being closed under the action of $\regring$ 
implies that for any $q \in Q_{I \rightarrow J}$ the map $\alpha_q$ is 
a zero map. Which by Proposition $\ref{prps-divisors-of-zeroes}$ implies 
$y \in B_{q}$. 

The support of each $M_\chi$ consists only of 
the prime divisors in $\mathfrak{E}$
($\cite{Logvinenko-Natural-G-Constellation-Families}$, Prop. 3.14 and
3.15). The same is true of the principal divisors $(x_i)$ for their 
images in $\mathbb{C}^n/G$ are the coordinate hyperplanes $x_i^{|G|} = 0$.
Therefore, by their defining equation $\eqref{eqn-divisor-of-zeroes}$, 
the support of each of the divisors of zeroes $B_q$ of $\mathcal{M}$ 
consists also only of the prime divisors in $\mathfrak{E}$. As $y$ lies 
on all $B_q$ with $q \in Q_{I \rightarrow J}$, $y$ must lie on at 
least one divisor in $\mathfrak{E}$. But, as $y \in U$, $y$ also lies 
on at most one divisor in $\mathfrak{E}$. Denote this
unique divisor by $E$, then  
\begin{align} \label{eqn-I-is-a-sub-R-module-condition}
q \in Q_{I \rightarrow J} \Rightarrow E \subset B_q.
\end{align}
Define a new $G$-Weil divisor set $\{M'_\chi\}$ by
setting $M'_\chi$ to be $M_\chi$ if $\chi \in I$ and 
$M_\chi + E$ if $\chi \in J$. Then divisors $\{ B'_q \}$ defined from 
$\{M'_\chi\}$ by equations $\eqref{eqn-divisor-of-zeroes}$ can be 
expressed as
\begin{align}\label{eqn-new-divisors-of-zeroes}
B'_q = 
\begin{cases}
B_q &\text{ if } q \in Q_{I \rightarrow I}, Q_{J \rightarrow J} \\
B_q + E &\text{ if } q \in Q_{J \rightarrow I} \\
B_q - E &\text{ if } q \in Q_{I \rightarrow J}
\end{cases}.
\end{align}
Since $\{B_q\}$ are all effective $\eqref{eqn-new-divisors-of-zeroes}$ 
and $\eqref{eqn-I-is-a-sub-R-module-condition}$ imply that $\{B'_q\}$
are also all effective. Therefore $\bigoplus \mathcal{L}(-M'_\chi)$
is a normalized $\gnat$-family. But  
\begin{align}
w_\theta(\mathcal{M}') =
w_\theta(\mathcal{M}) + \sum_{\chi \in J}
\theta(\chi) 
\end{align}
which contradicts the maximality of $w_\theta(\mathcal{M})$ since 
$\sum_{\chi \in J} \theta(\chi) = - \theta(A') > 0$.  

For the second claim let $\mathcal{N} = \bigoplus \mathcal{L}(-N_\chi)$
be another normalized family $\theta$-semistable over $U$. Let $B'_q$ 
be divisors of zeroes of $\mathcal{N}$. Then 
\begin{align} \label{eqn-divisors-of-zero-difference}
B_q - B'_q = (M_{tq} - N_{tq}) - (M_{hq} - N_{hq}).
\end{align}

Take any $E' \in \mathfrak{E}$ such that the sets 
$\{m_{\chi, E'}\}$ and $\{n_{\chi, E'}\}$ of the coefficients of $E'$ 
in $\{M_\chi\}$ and $\{N_\chi\}$ are distinct. Then 
$J' = \{\chi \in G^\vee \;|\; n_{\chi,E'} > m_{\chi,E'} \}$ is 
a non-empty proper subset of $G^\vee$. Denote by $I'$ its complement. 
For any $q \in Q_{I' \rightarrow J'}$ the coefficient of $E'$
in the RHS of $\eqref{eqn-divisors-of-zero-difference}$ is strictly
positive. As $B'_q$ is effective we conclude that 
$q \in Q_{I' \rightarrow J'}$ implies  $E' \subset B_q$.
So for any $y \in E'$ the restriction $(\bigoplus_{\chi \in I'}
\mathcal{L}(M_\chi))|_y$ is a sub-$\twalg$-module of
$\mathcal{M}_{|y}$. But as $\mathcal{M}$ is $\theta$-semistable on $U$
and as $U \cap E' \neq \emptyset$
we must have $\sum_{\chi \in I'} \theta(\chi) \geq 0$. Similarly 
if $q \in Q_{J' \rightarrow I'}$, then the RHS of
$\eqref{eqn-divisors-of-zero-difference}$ is strictly negative, so 
$E' \subset B'_q$ and $\theta$-semistability of $\mathcal{N}$ implies
$\sum_{\chi \in J'} \theta(\chi) = - \sum_{\chi \in I'} \theta(\chi)
\geq 0$. Therefore $\sum_{\chi \in I'} \theta(\chi) = 0$ and $\theta$
is not generic. 
\end{proof}

The fine moduli space $M_\theta$ of $\theta$-stable 
$G$-constellations can be constructed via GIT theory, together with 
the universal family $\mathcal{M}_\theta$. The Hilbert-Chow morphism
$\pi_\theta$ of $\mathcal{M}_\theta$ is projective. As the universal
family is defined up to an equivalence of families, that is up to 
a twist by a line bundle, we can assume $\mathcal{M}_\theta$ 
to be normalised.

Assume for the rest of this section that $n=3$. If $\theta$ is 
generic, then $M_\theta$ is a projective crepant resolution 
of $\mathbb{C}^3/G$ and $\mathcal{M}_\theta$ is 
everywhere orthogonal in all degrees. As any two crepant resolutions
of a canonical treefold are connected by a chain of flops, $M_\theta$
and $Y$ are isomorphic outside of a codimension $2$ subset. The
maps $Y \xrightarrow{\pi} \mathbb{C}^3/G$ and $M_\theta
\xrightarrow{\pi_\theta} \mathbb{C}^3/G$ fix a choice of a birational
isomorphism between $Y$ and $M_\theta$. This, as described in Section 
$\ref{section-direct-transforms-and-stable-families}$, defines
a notion of direct transforms between $Y$ and $M_\theta$.  

\begin{corollary} \label{cor-direct-transform-of-stable-formula}
Let $\theta \in \Theta$ be generic. Let $\mathcal{M}$ be
the unique normalized $\gnat$-family on $Y$ which maximizes 
the map $w_\theta$. Then $\mathcal{M}$ is isomorphic to 
the direct transform of $\mathcal{M}_\theta$ from $M_\theta$ to $Y$.
\end{corollary}

\begin{proof}
By the first claim of Proposition $\ref{prps-stability-criterion}$,
$\mathcal{M}$ is $\theta$-stable on $U$. So, by its definition, is 
the direct transform of $\mathcal{M}_{\theta}$ to $Y$. Hence, by 
the second claim of Proposition $\ref{prps-stability-criterion}$, 
$\mathcal{M}$ and the direct transform of $\mathcal{M}_{\theta}$ must 
be isomorphic.
\end{proof}

\section{Non-projective example} \label{section-non-projective-example}

In this section we give an application of the 
Theorem \ref{theorem-main-theorem} whereby we 
construct explicitly a derived McKay correspondence  
for a choice of an abelian $G \subset \gsl_3(\mathbb{C})$ and of a 
non-projective crepant resolution $Y$ of $\mathbb{C}^3/G$. 

\subsection{The group}
\label{subsection-the-group}

We set the group $G$ to be $\frac{1}{6}(1,1,4) \oplus \frac{1}{2}(1,0,1)$. 
That is, the image in $\gsl_3(\mathbb{C})$ of the product 
$\mu_6 \times \mu_2$ of groups of $6$th and $2$nd roots of unity, 
respectively, under the embedding:
\begin{align} \label{eqn-G-defn}
(\xi_1, \xi_2) \mapsto 
\begin{pmatrix}
\xi_1 \xi_2 & & \\
& \xi_1 & \\
& & \xi_1^{4} \xi_2
\end{pmatrix}.
\end{align}

We denote by $\chi_{i,j}$ the character of $G$ induced by 
$(\xi_1, \xi_2) \mapsto \xi_1^i \xi_2^j$. 

Calculating the McKay quiver of $G$ (cf. Section 
\ref{section-the-mckay-quiver-of-G}), we obtain:
\begin{center}
\includegraphics[scale=0.15]{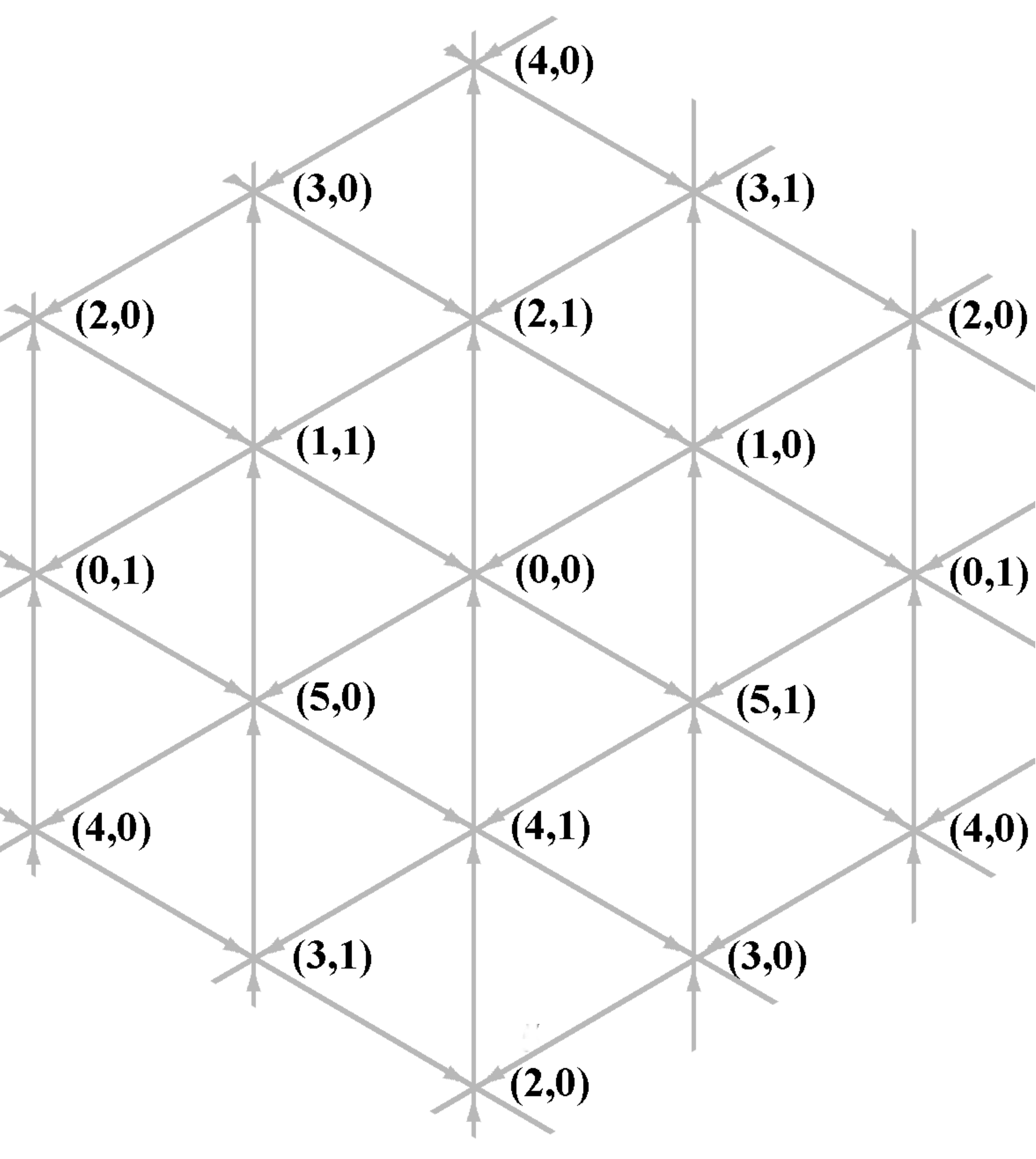}
\end{center}
The way we've chosen to depict the McKay quiver reflects 
the fact that it has a universal cover quiver naturally embedded
into $\mathbb{R}^2$. This point of view 
will not be essential for our argument but a curious reader should 
consult $\cite{Craw-Ishii-02}$, Section 10.2 and 
$\cite{Logvinenko-thesis}$, Section 6.4.

\subsection{The resolution}
\label{subsection-toric-generalities}

We define the crepant resolution $Y$ of
$\mathbb{C}^3 / G$ using methods of toric geometry.
For the specifics related to $G$-constellations see $\cite{Logvinenko-Families-of-G-constellations-over-resolutions-of-quotient-singularities}$, Section 3.  

We define the relevant notation. The embedding $\eqref{eqn-G-defn}$
defines a surjection of torii 
\begin{align}\label{seq-g-torus}
\xymatrix{
0 \ar[r] &
G \ar[r] &
(\mathbb{C}^*)^3 \ar[r] &
T \ar[r] &
0
}.
\end{align}
Applying $\homm(\bullet, \mathbb{C}^*)$ to $\eqref{seq-g-torus}$ 
we obtain the character lattices of the torii:
\begin{align} \label{seq-monom-char}
\xymatrix{
0 \ar[r] &
M \ar[r] &
\mathbb{Z}^3 \ar[r]^{\rho} &
G^\vee \ar[r] &
0
}.
\end{align}
Given any character $m = (k_1, k_2, k_3) \in \mathbb{Z}^3$ of 
$(\mathbb{C}^*)^3$ we denote by $x^{m}$ the Laurent monomial 
$x_1^{k_1} x_2^{k_2} x_3^{k_3}$ in $R$. Applying 
$\homm(\bullet, \mathbb{Z})$ to $\eqref{seq-monom-char}$ 
we obtain the dual lattices
\begin{align*}
\xymatrix{
0 \ar[r] &
(\mathbb{Z}^3)^\vee \ar[r] &
N \ar[r]&
\ext^1(G^\vee, \mathbb{Z}) \ar[r] &
0
}.
\end{align*}
Let $e_1, e_2, e_3$ be the basis of $(\mathbb{Z}^3)^\vee$ dual to
$x_1, x_2, x_3$. The dual lattice $N$ is generated over 
$(\mathbb{Z}^3)^\vee$ by $\frac{1}{6}(1,1,4)$ and
$\frac{1}{2}(1,0,1)$. The quotient space $\mathbb{C}^3/G$ is the toric
variety given by a single cone $\sigma_{\geq 0} = \sum
\mathbb{R}_{\geq 0} e_i$ in $N$. Let $Y$ be the toric variety whose 
fan $\mathfrak{F}$ in $N$ is the subdivision of $\sigma_{\geq 0}$ which 
triangulates the junior simplex 
$\Delta = \{ (k_1, k_2, k_3) \in \sigma_{\geq 0}
\; |\; \sum k_i = 1 \}$ as depicted below
\begin{center}
\includegraphics[scale=0.20]{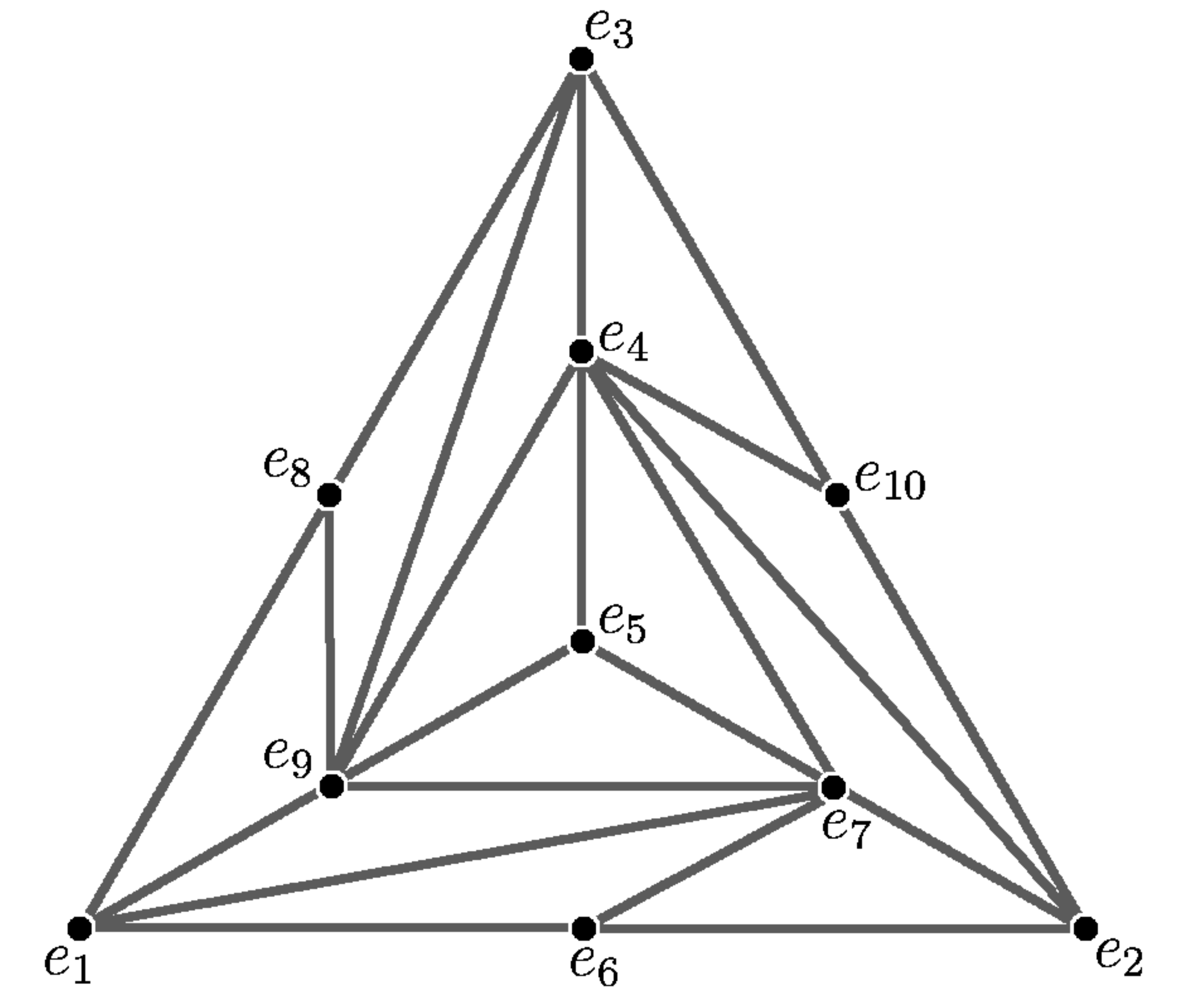}
\end{center}
where by $e_i$ we denote the following elements of $N$ 
\begin{align} \label{eqn-e_i}
\begin{tabular}{l l l}
$e_1 = (1,0,0)$  &  $e_2 = (0,1,0)$  &  $e_3 = (0,0,1)$ \\
$e_4 = \frac{1}{6}(1,1,4)$  &  $e_5 = \frac{1}{3}(1,1,1)$ 
& $e_6 = \frac{1}{2}(1,1,0)$ \\
$e_7 = \frac{1}{6}(1,4,1)$  &  $e_8 = \frac{1}{2}(1,0,1)$ 
& $e_9 = \frac{1}{6}(4,1,1)$ \\
$e_{10} = \frac{1}{2}(0,1,1).$  &  &
\end{tabular}
\end{align}
Denote by $\pi$ the map $Y \rightarrow \mathbb{C}^3/G$ defined 
by the inclusion of $\mathfrak{F}$ into $\sigma_{\geq 0}$. All  
the maximal cones of $\mathfrak{F}$ are basic in $N$, so 
$Y$ is smooth. The generators $e_i$ of the rays of $\mathfrak{F}$
lie in $\Delta$, so the map $\pi$ is crepant(\cite{YPG87}, Prop.
4.8). Finally, the argument of
\cite{KKMSD-ToroidalEmbeddingsI}, Chapter III, \S 2E, Example 2 
shows that $\pi$ is non-projective. 

The quotient torus $T$ acts on $Y$ and to each $k$-dimensional cone 
$\sigma$ in $\mathfrak{F}$ corresponds 
a $(3-k)$-dimensional orbit of $T$. We denote it by $S_\sigma$
and denote by $E_\sigma$ the closure of $S_\sigma$, it is
the union of all orbits $S_{\sigma'}$ with $\sigma \subseteq \sigma'$.
For each cone $\left<e_i\right>$ in the fan $\mathfrak{F}$, we denote 
by $S_{i}$ the codimension $1$ orbit $S_{\left<e_i\right>}$ and by $E_{i}$ 
the divisor $E_{\left<e_i\right>}$. Similarly we use $S_{i,j}$ and
$E_{i,j}$ for the codimension $2$ orbit $S_{\left<e_i,e_j\right>}$ and 
the surface $E_{\left<e_i,e_j\right>}$ and we use $E_{i,j,k}$ for 
the toric fixed point $E_{\left<e_i, e_j, e_k\right>}$. 

\subsection{The family}
\label{section-the-family}

The map $Y \xrightarrow{\pi} \mathbb{C}^3/G$ defines the notion
of $G$-Weil divisors on $Y$. Any 
normalized $\gnat$-family on $Y \xrightarrow{\pi} \mathbb{C}^3/G$ 
is of the form $\bigoplus_{\chi \in G^\vee} \mathcal{L}(-D_\chi)$ for 
some $G$-Weil divisors $D_\chi$ with $D_{\chi_{0,0}} = 0$. Moreover, as explained 
in \cite{Logvinenko-Natural-G-Constellation-Families}, Section 3.5, 
there exists the \it maximal shift family \rm $\oplus \mathcal{L}(- M_\chi)$
such that for any other normalized $\gnat$-family $\oplus \mathcal{L}(- D_\chi)$
we have 
\begin{align}\label{eqn-maxshift-fam}
M_\chi \geq D_\chi 
\end{align}
for all $\chi \in G^\vee$. We denote this family by $\mathcal{F}$ and 
shall prove it to satisfy the assumptions of Corollary 
\ref{cor-3d-application-of-main-theorem}.

In the notation of Section \ref{subsection-toric-generalities}
each divisor $M_\chi$ is of form $\sum q_{\chi,i} E_i$.
The coefficients $q_{\chi,i}$ can be calculated via formula
\begin{align} \label{eqn-maximal-shift-set}
	q_{\chi,i} = 
	\inf \{ e_i(m) \; | \; m \in \sigma_{\geq 0}^\vee 
\cap \rho^{-1}(\chi) \}.
\end{align}
A detailed example of such calculation can be seen in 
$\cite{Logvinenko-Families-of-G-constellations-over-resolutions-of-quotient-singularities}$,
Example 4.21. In our case, we obtain $q_{\chi,i}$ to be: 
\begin{align} \label{eqn-maximal-shift-family}
\text{
\begin{scriptsize}
\begin{tabular}{|r|r|r|r|r|r|r|r||r|r|r|r|r|r|r|r|r|r|}
\hline
$\chi \setminus i $ & $ 4 $ & $ 5 $ & $ 6 $ & $ 7 $ & $ 8 $ & $ 9 $ & $ 10$ &
$\chi \setminus i $ & $ 4 $ & $ 5 $ & $ 6 $ & $ 7 $ & $ 8 $ & $ 9 $ & $ 10$ \\ 
\hline
& & & & & & & & & & & & & & & \\
$\chi_{0,0}$ & $ 0 $ & $ 0 $ & $ 0 $ & $ 0 $ & $ 0 $ & $ 0 $ & $ 0
$ &
$\chi_{2,0}$ & $ \frac{2}{6} $ & $ \frac{4}{6} $ & $ 0 $
& $ \frac{2}{6} $ & $ 0 $ & $ \frac{2}{6} $ & $ 0 $ \\ 
& & & & & & & & & & & & & & & \\
$\chi_{4,0}$ & $ \frac{4}{6} $ & $ 1\frac{2}{6} $ & $ 0 $
& $ \frac{4}{6} $ & $ 0 $ & $ \frac{4}{6} $ & $ 0 $ & 
$\chi_{1,1}$ & $ \frac{1}{6} $ & $ \frac{2}{6} $ & $ \frac{3}{6} $
& $ \frac{1}{6} $ & $ \frac{3}{6} $ & $ \frac{4}{6} $ & $ 0 $ \\ 
& & & & & & & & & & & & & & & \\
$\chi_{1,0}$ & $ \frac{1}{6} $ & $ \frac{2}{6} $ & $ \frac{3}{6} $
& $ \frac{4}{6} $ & $ 0 $ & $ \frac{1}{6} $ & $ \frac{3}{6} $ & 
$\chi_{4,1}$ & $ \frac{4}{6} $ & $ \frac{2}{6} $ & $ 0 $
& $ \frac{1}{6} $ & $ \frac{3}{6} $ & $ \frac{1}{6} $ & $ \frac{3}{6} $ \\ 
& & & & & & & & & & & & & & & \\
$\chi_{3,1}$ & $ \frac{3}{6} $ & $ 1 $ & $ \frac{3}{6} $
& $ \frac{3}{6} $ & $ \frac{3}{6} $ & $ 1 $ & $ 0 $ & 
$\chi_{3,0}$ & $ \frac{3}{6} $ & $ 1 $ & $ \frac{3}{6} $
& $ 1 $ & $ 0 $ & $ \frac{3}{6} $ & $ \frac{3}{6} $ \\ 
& & & & & & & & & & & & & & & \\
$\chi_{0,1}$ & $ 1 $ & $ 1 $ & $ 0 $
& $ \frac{3}{6} $ & $ \frac{3}{6} $ & $ \frac{3}{6} $ & $ \frac{3}{6} $ & 
$\chi_{5,1}$ & $ \frac{5}{6} $ & $ \frac{4}{6} $ & $ \frac{3}{6} $
& $ \frac{5}{6} $ & $ \frac{3}{6} $ & $ \frac{2}{6} $ & $ 0 $ \\ 
& & & & & & & & & & & & & & & \\
$\chi_{5,0}$ & $ \frac{5}{6} $ & $ \frac{4}{6} $ & $ \frac{3}{6} $
& $ \frac{2}{6} $ & $ 0 $ & $ \frac{5}{6} $ & $ \frac{3}{6} $ & 
$\chi_{2,1}$ & $ \frac{2}{6} $ & $ \frac{4}{6} $ & $ 0 $
& $ \frac{5}{6} $ & $ \frac{3}{6} $ & $ \frac{5}{6} $ & $ \frac{3}{6} $ \\ 
& & & & & & & & & & & & & & & \\
\hline
\end{tabular}
\end{scriptsize}
}
\end{align}

The principal $G$-Weil divisors $(x_k)$ can 
be calculated with a formula 
\begin{align}
(x_i) = \frac{1}{12} \sum^{10}_{j=1} e_j(x_i^{12}) E_j, 
\end{align}
cf. \cite{Logvinenko-Families-of-G-constellations-over-resolutions-of-quotient-singularities},
Prop. 3.2. In our case we obtain:
\begin{align} \label{eqn-principal-divisors}
\text{
\begin{small}
\begin{tabular}{l l l}
& $(x_1) = E_1 + \frac{1}{6}E_4 + \frac{1}{3}E_5 + \frac{1}{2}E_6 +  
\frac{1}{6}E_7 + \frac{1}{2}E_8 + \frac{4}{6}E_9$ &\\
& & \\
& $(x_2) = E_2 + \frac{1}{6}E_4 + \frac{1}{3}E_5 + \frac{1}{2}E_6 +  
\frac{4}{6}E_7 + \frac{1}{6}E_9 + \frac{1}{2}E_{10}$ &\\
& & \\  
& $(x_3) = E_3 + \frac{4}{6}E_4 + \frac{1}{3}E_5 + \frac{1}{6}E_7 +
\frac{1}{2}E_8 + \frac{1}{6}E_9 + \frac{1}{2}E_{10}$ & 
\end{tabular}
\end{small}
}
\end{align}
Substituting the data of $\eqref{eqn-principal-divisors}$ and
$\eqref{eqn-maximal-shift-family}$ into the formula
$\eqref{eqn-divisor-of-zeroes}$ we calculate for every arrow 
of the McKay quiver its divisor of zeroes in $\mathcal{F}$:
\begin{align} \label{eqn-divisors-of-zeroes}
\text{
\begin{footnotesize}
\begin{tabular}{l l}
$ B_{\chi_{0,0}, 1} = E_1 $ &
$ B_{\chi_{1,1}, 1} = E_1 + E_4 + E_5 + E_6 + E_7 + E_8 + E_9 $ \\
$ B_{\chi_{0,0}, 2} = E_2 $ & 
$ B_{\chi_{1,1}, 2} = E_2 + E_6 + E_7 $ \\
$ B_{\chi_{0,0}, 3} = E_3 $ & 
$ B_{\chi_{1,1}, 3} = E_3 + E_4 + E_8 $ \\
$ B_{\chi_{4,0}, 1} = E_1 $ & 
$ B_{\chi_{1,0}, 1} = E_1 + E_6 + E_9$ \\
$ B_{\chi_{4,0}, 2} = E_2 $ & 
$ B_{\chi_{1,0}, 2} = E_2 + E_4 + E_5 + E_6 + E_7 + E_9 + E_{10} $ \\
$ B_{\chi_{4,0}, 3} = E_3 $ & 
$ B_{\chi_{1,0}, 3} = E_3 + E_4 + E_{10} $ \\
$ B_{\chi_{2,0}, 1} = E_1 + E_5 + E_9 $ & 
$ B_{\chi_{4,1}, 1} = E_1 + E_8 + E_9 $ \\
$ B_{\chi_{2,0}, 2} = E_2 + E_5 + E_7 $ & 
$ B_{\chi_{4,1}, 2} = E_2 + E_7 + E_{10} $ \\
$ B_{\chi_{2,0}, 3} = E_3 + E_4 + E_5 $ & 
$ B_{\chi_{4,1}, 3} = E_3 + E_4 + E_5 + E_7 + E_8 + E_9 + E_{10} $ \\
$ B_{\chi_{5,1}, 1} = E_1 + E_6 + E_8 + E_9 $ &
$ B_{\chi_{3,1}, 1} = E_1 + E_6 + E_8 + E_9  $ \\
$ B_{\chi_{5,1}, 2} = E_2 + E_6 $ & 
$ B_{\chi_{3,1}, 2} = E_2 + E_5 + E_6 + E_7 + E_9 $ \\
$ B_{\chi_{5,1}, 3} = E_3 + E_8 $ & 
$ B_{\chi_{3,1}, 3} = E_3 + E_4 + E_5 + E_8 + E_9 $ \\
$ B_{\chi_{5,0}, 1} = E_1 + E_6  $ & 
$ B_{\chi_{3,0}, 1} = E_1 + E_5 + E_6 + E_7 + E_9  $ \\
$ B_{\chi_{5,0}, 2} = E_2 + E_6 + E_7 + E_{10} $ & 
$ B_{\chi_{3,0}, 2} = E_2 + E_6 + E_7 + E_{10} $ \\
$ B_{\chi_{5,0}, 3} = E_3 + E_{10} $ & 
$ B_{\chi_{3,0}, 3} = E_3 + E_4 + E_5 + E_7 + E_{10} $ \\
$ B_{\chi_{2,1}, 1} = E_1 + E_8  $ & 
$ B_{\chi_{0,1}, 1} = E_1 + E_4 + E_5 + E_8 + E_9 $ \\
$ B_{\chi_{2,1}, 2} = E_2 + E_{10} $ & 
$ B_{\chi_{0,1}, 2} = E_2 + E_4 + E_5 + E_7 + E_{10}  $ \\
$ B_{\chi_{2,1}, 3} = E_3 + E_4 + E_{8} + E_{10} $ & 
$ B_{\chi_{0,1}, 3} = E_3 + E_4 + E_8 + E_{10}.$ \\
\end{tabular}
\end{footnotesize}
}
\end{align}

\subsection{A sample calculation} \label{subsection-sample-calculation}

Corollary $\ref{cor-orthogonality-criterion}$ together with 
the table $\eqref{eqn-divisors-of-zeroes}$ 
are all that we need to check any two $G$-constellations in 
$\mathcal{F}$ for the degree $0$ orthogonality.
Below we give an example of a calculation which verifies 
that any point on the torus orbit $S_8$ and any 
point on the torus orbit $S_{1,7}$ are orthogonal in degree $0$ 
in $\mathcal{F}$.  

Let $a$ be any point of $S_8$. Then $a$ lies on no divisor $E_i$ 
other than $E_8$. Hence $a \in B_q$
if and only if $E_8 \subset B_q$. Let $A$ be the fiber 
of $\mathcal{F}$ at $a$ and $\{\alpha_q\}$ be the corresponding 
representation of the McKay quiver. By Proposition 
$\ref{prps-divisors-of-zeroes}$ for any arrow $q$ the map $\alpha_q$
is a zero map if and only if $E_8 \in B_q$. On Figure 4
we use the table $\eqref{eqn-divisors-of-zeroes}$ and mark all
the zero-maps in $\{\alpha_q\}$ by
drawing a line through the corresponding arrow of the McKay quiver.
Similarly if $b$ is a point of $S_{1,7}$ then $b$ lies on no 
$E_i$ other than $E_1$ and $E_7$. Let $B$ be the fiber
of $\mathcal{F}$ at $b$ and $\{\beta_q\}$ be the corresponding
representation. As above $\beta_q$
is a zero-map if and only if either $E_1$ or $E_7$ belongs 
to $B_q$. On Figure 5 we mark all the zero-maps $\{\beta_q\}$.
 
\begin{center}
\begin{tabular}{c c}
\includegraphics[scale=0.130]{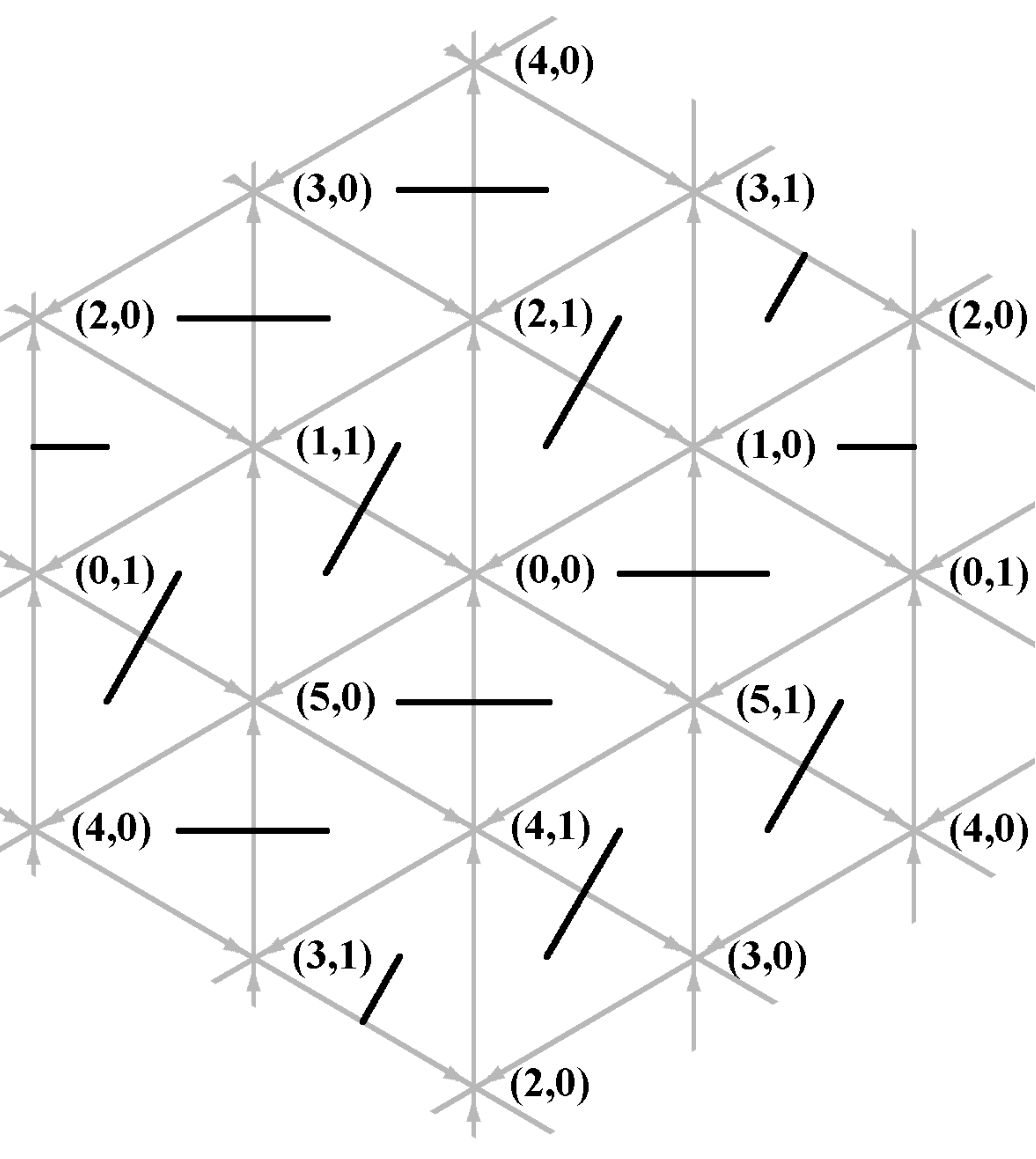} &
\includegraphics[scale=0.130]{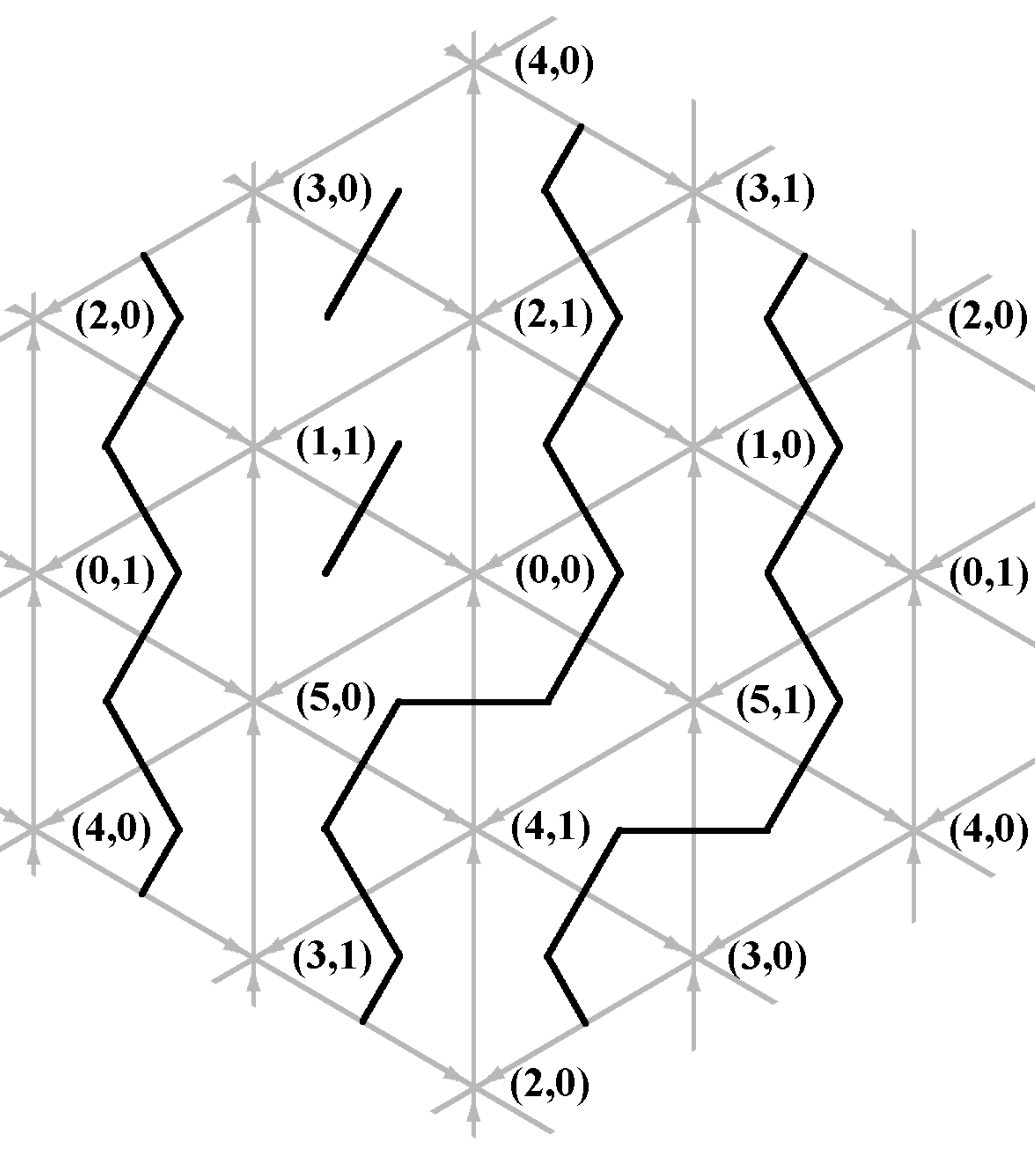} \\
\bf Figure 4 \rm & \bf Figure 5 \rm
\end{tabular}
\end{center}

\begin{center}
\begin{tabular}{c c}
\includegraphics[scale=0.130]{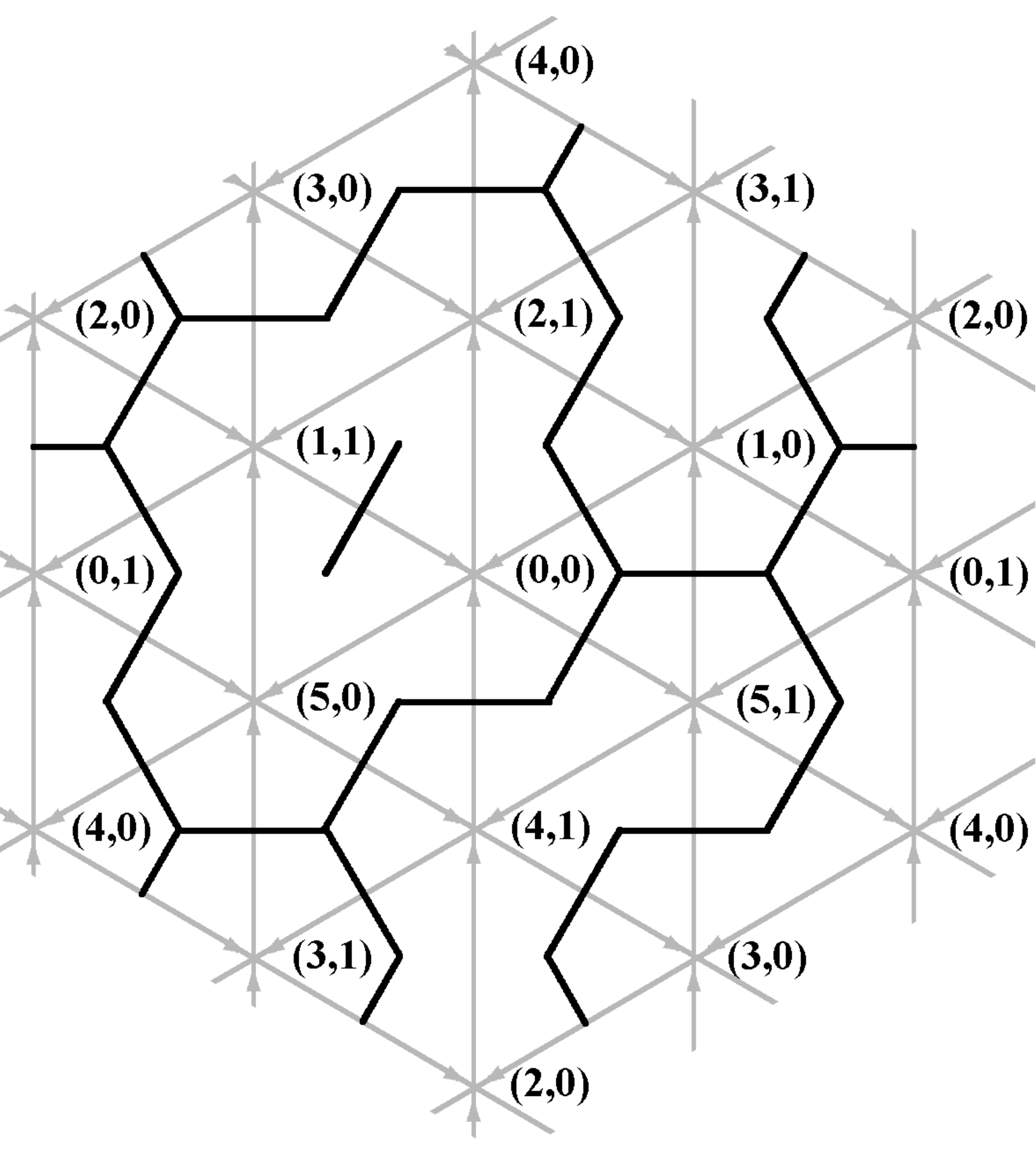} & 
\includegraphics[scale=0.130]{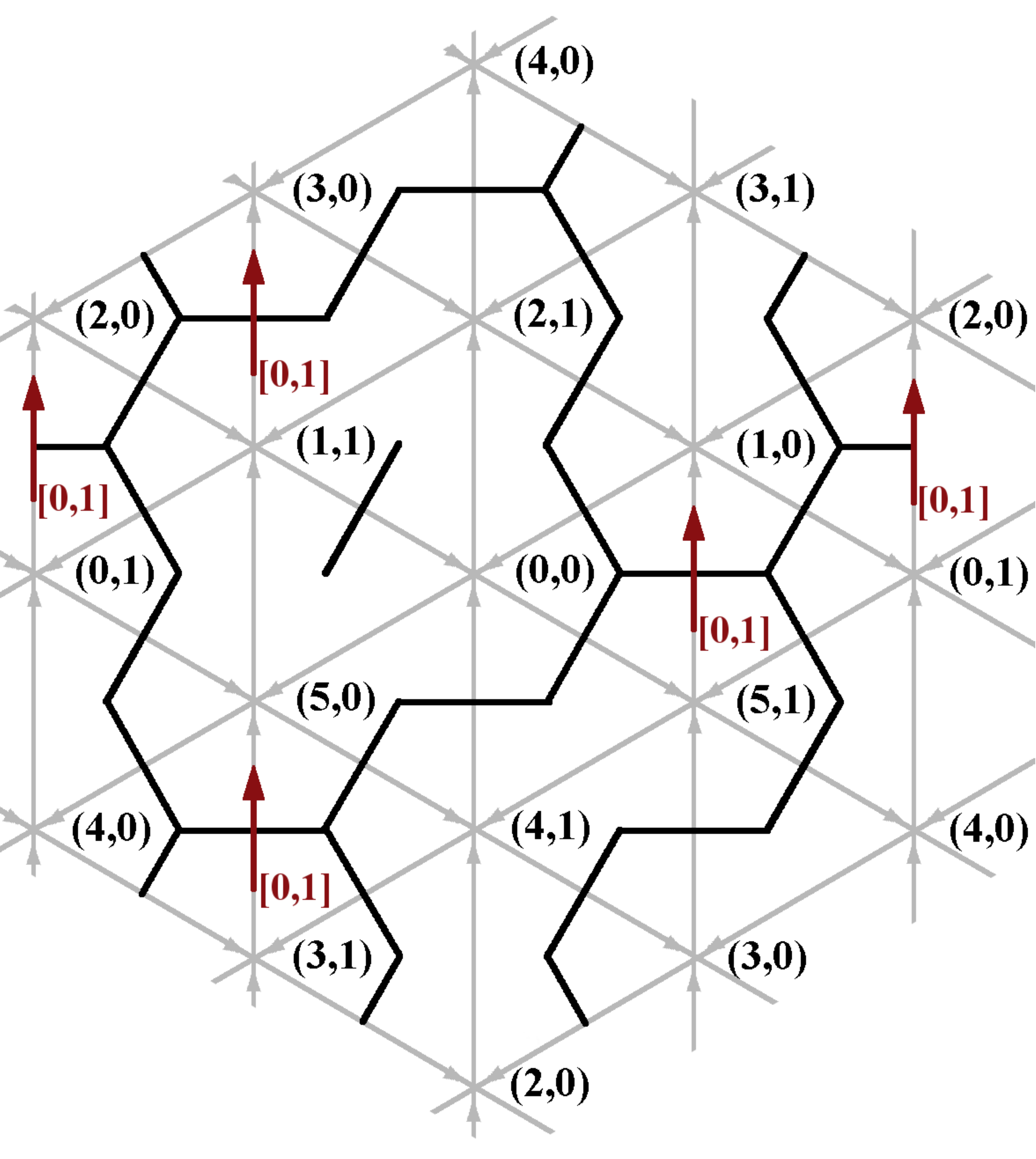} \\
\bf Figure 6 \rm &
\bf Figure 7 \rm
\end{tabular}
\end{center}

On Figure $6$ we combine the markings of Figures $4$ and $5$.
The arrows left unmarked are the arrows of type $[1,1]$ with respect
to the pair $A,B$ (Def.~$\ref{defn-arrow-type}$).
It is clear that the components path-connected by 
$[1,1]$-arrows are: $\{\chi_{0,0},
\chi_{2,1}, \chi_{5,0}, \chi_{1,1} \}$, $\{\chi_{5,1}, \chi_{4,1}, 
\chi_{2,0} \}$, $\{ \chi_{1,0}, \chi_{3,1} \}$ and $\{ \chi_{0,1},
\chi_{4,0}, \chi_{3,0}\}$. Now, with Cor. 
$\ref{cor-orthogonality-criterion}$ in mind, we search the borders of these
four regions for the $[1,0]$ and $[0,1]$-arrows. The $[1,0]$-arrows 
are the ones unmarked on Figure $4$ but marked 
on Figure $5$ and vice versa for $[0,1]$. On Figure $7$ we've 
marked on the border of each region an incoming and 
an outgoing $[0,1]$-arrow. By Cor.
$\ref{cor-orthogonality-criterion}$ we see that $A$ 
and $B$ are orthogonal in degree $0$. 

\subsection{Final calculations}

We now claim that $\mathcal{F}$ is the direct transform of the
universal family of $G$-clusters on $G$-$\hilb(\mathbb{C}^3)$. 
In the notation of Section $\ref{subsection-theta-stability}$ 
define $\theta_+ \in \Theta$ by $\theta_+(\chi_{0,0}) = 1 - |G|$ and 
$\theta_+(\chi) = 1$ for $\chi \neq \chi_{0,0}$. Evidently $\theta_+$ 
is generic. It follows from the original observation by Ito and 
Nakajima in $\cite{ItoNakajima98}$, \S3, that $G$-clusters can be identified 
with $\theta_+$-stable $G$-constellations, thus identifying
$G$-$\hilb(\mathbb{C}^3)$ with the fine moduli space $M_{\theta_+}$.
On the other hand, inequalities \eqref{eqn-maxshift-fam} imply  
that $\mathcal{F}$ maximizes $\omega_{\theta_+}$ on 
$Y \xrightarrow{\pi} \mathbb{C}^3/G$. Hence, by Corollary 
\ref{cor-direct-transform-of-stable-formula}, $\mathcal{F}$ is
the direct transform of $\mathcal{M}_{\theta_+}$ from 
$G$-$\hilb(\mathbb{C}^3)$ to $Y$. 

For a detailed description of an algorithm which allows one 
to calculate the toric fan of $G$-$\hilb(\mathbb{C}^3)$ 
see in $\cite{Craw02}$. For our group $G$ we obtain:
\begin{center}
\includegraphics[scale=0.18]{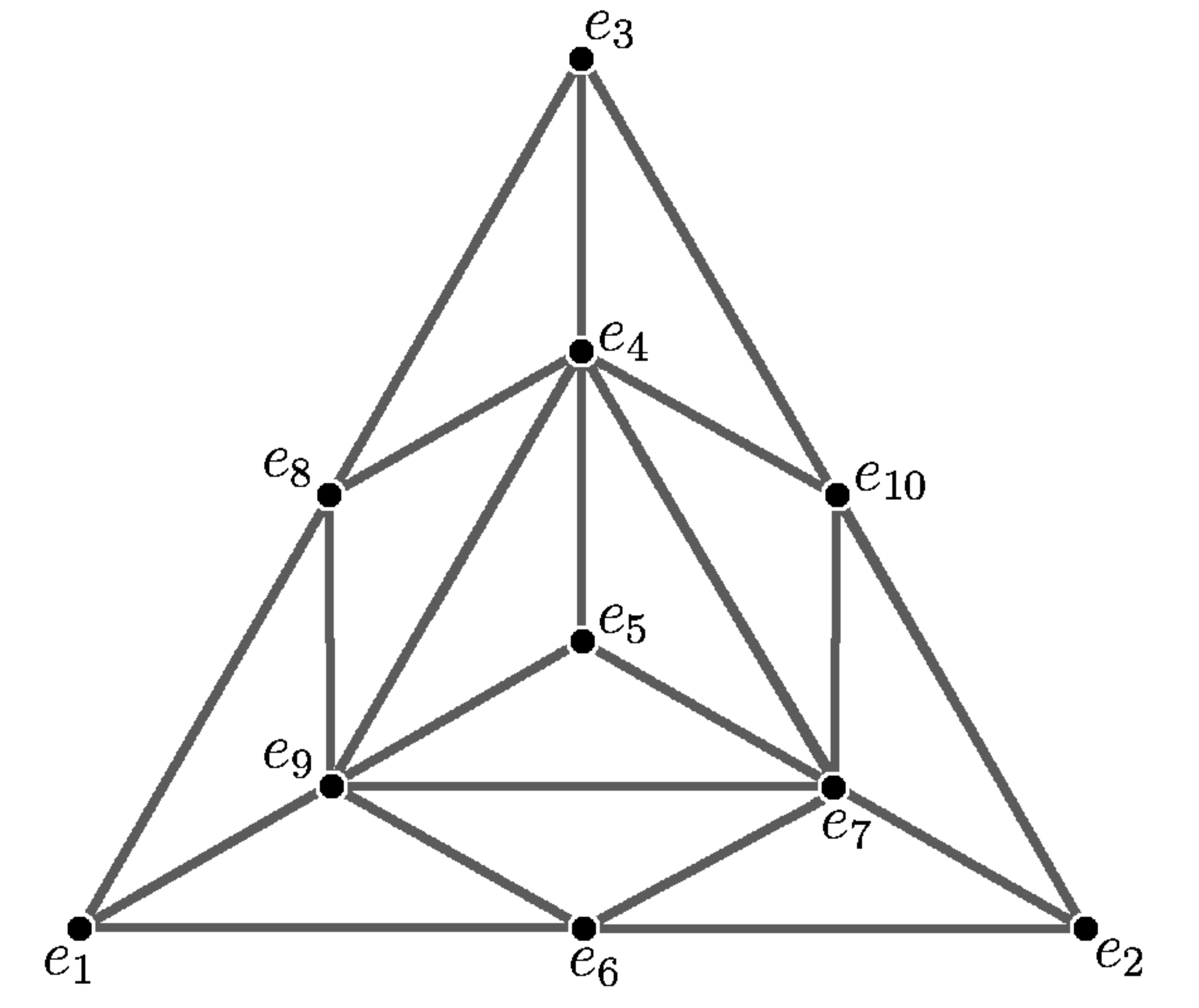}
\end{center}
\begin{center}
\bf Figure 8 \rm
\end{center}

The general points of an exceptional surface $E_i$, as per the
statement of Corollary \ref{cor-3d-application-of-main-theorem}, are
precisely the codimension $1$ torus orbit $S_i$. Similarly, the
general points of an exceptional curve $E_i \cap E_j$ are precisely
the codimension $2$ torus orbit $S_{i,j}$. Comparing Figure $8$ with 
the fan of $Y$ on Figure $3$ we see that the only codimension
$1$ or $2$ torus orbits in $Y$ whose corresponding cones 
aren't also contained in the fan of $G$-$\hilb(\mathbb{C}^3)$ are 
$S_{1,7}$, $S_{2,4}$ and $S_{3,9}$. The argument in Section 
\ref{section-direct-transforms-and-stable-families}
reduces verifying that $\mathcal{F}$ satisfies the conditions 
of Corollary \ref{cor-3d-application-of-main-theorem},
to checking that each of these three orbits is orthogonal 
in degree $0$ in $\mathcal{F}$ to every codimension $1$
orbit $S_i$.

We claim that, in fact, it suffices to check it for just one of these
orbits. Let $\phi$ be the rotation of the fan of $Y$ around the ray $e_5$ 
which rotates Figure $2$ clockwise by $2\pi/3$. Let $\psi$ be the
rotation of the plane containing the McKay quiver on the Figure $3$ 
anti-clockwise by $2\pi/3$ with center at $\chi_{0,0}$. 
Observe that the permutation of the divisors $E_i$ defined by $\phi$
and the permutation of the arrows of the McKay quiver defined by
$\psi$ leave the numerical data 
$\eqref{eqn-divisors-of-zeroes}$ of divisors of zeroes 
of $\mathcal{F}$ invariant\footnote{This invariance is a consequence 
of the fan of $Y$
being symmetric and of $\mathcal{F}$ being intrinsically defined as 
the maximal shift family.}. It follows that the orthogonality
calculation of Section $\ref{subsection-sample-calculation}$ for any
pair of torus orbits $S,S'$ and the same calculation for 
$\phi(S), \phi(S')$ differ on Figures $4$-$7$ only by a rotation by
$\psi$. The claim now follows as the cones of $S_{1,7}$, $S_{2,4}$ and
$S_{3,9}$ are permuted by $\phi$. 

We choose to treat $S_{1,7}$. We repeat the calculation
of Section $\ref{subsection-sample-calculation}$
for $S_{1,7}$ and every other orbit $S_i$ and list below
the analogues of Figure $7$. From them, as elaborated in Section 
$\ref{subsection-sample-calculation}$, the reader could readily
ascertain the orthogonality in $\mathcal{F}$ of 
the torus orbits involved. 

We conclude, by Corollary \ref{cor-3d-application-of-main-theorem}, 
that the integral transform $\Phi_{\mathcal{F}}(- \otimes \rho_0)$ 
is an equivalence of categories $D(Y) \rightarrow D^G(\mathbb{C}^3)$ and 
that \it a posteriori \rm the family $\mathcal{F}$ is everywhere orthogonal 
in all degrees. 

\begin{center}
\begin{tabular}{c c}
\includegraphics[scale=0.135]{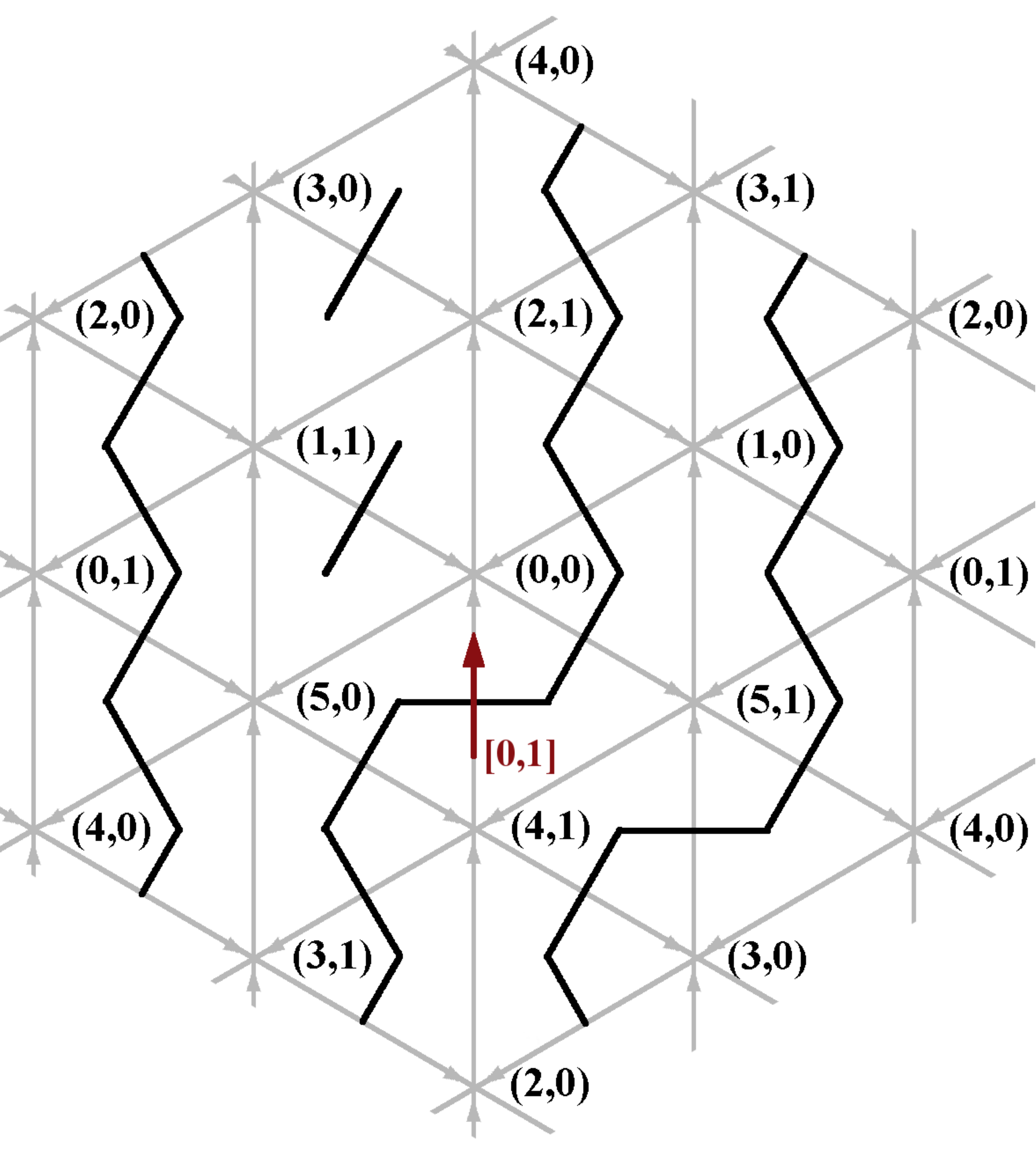} &
\includegraphics[scale=0.135]{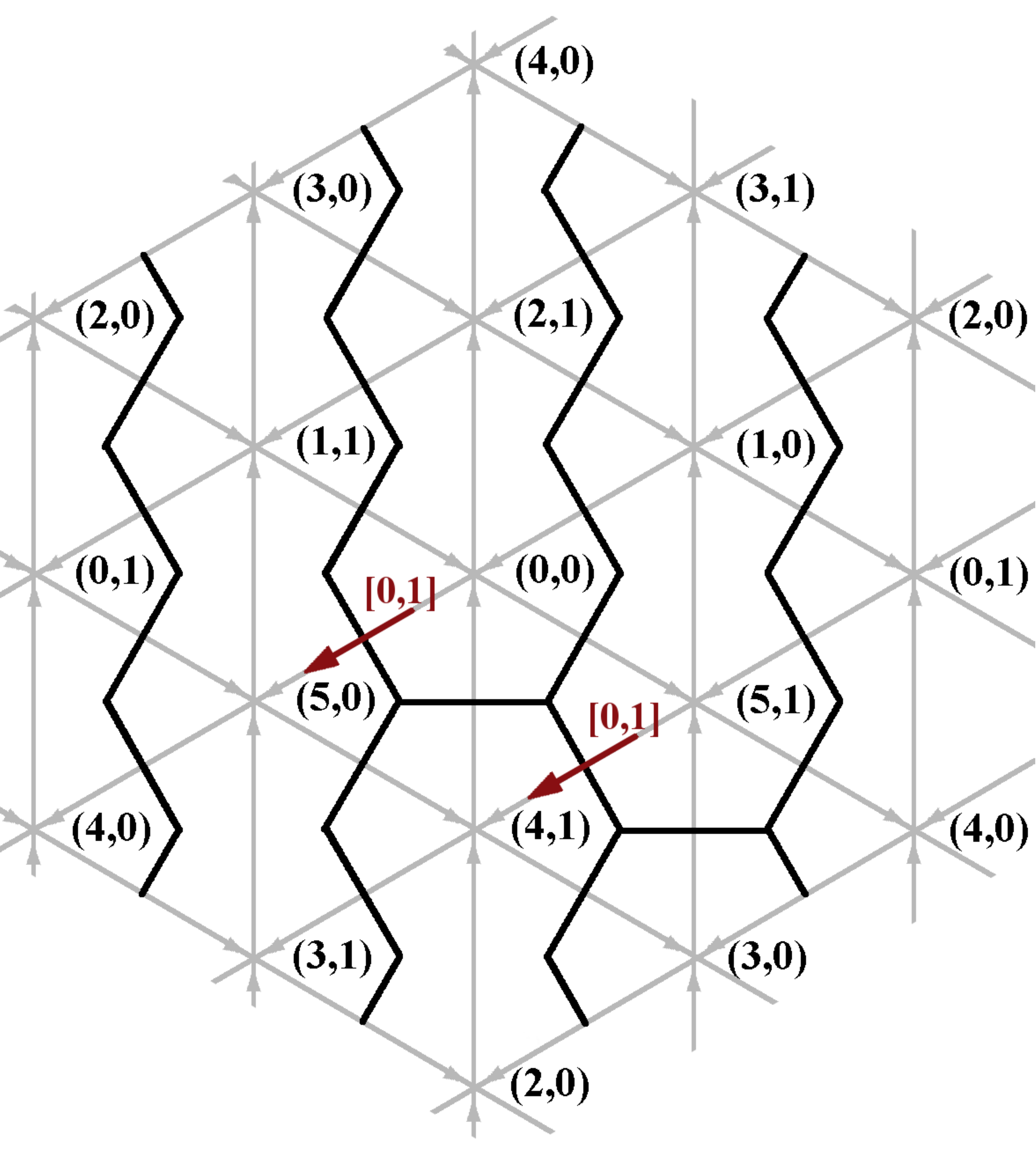} \\
$(S_{1},S_{1,7})$ and $(S_{7},S_{1,7})$ & $(S_{2},S_{1,7})$ 
\end{tabular}
\end{center}
\begin{center}
\begin{tabular}{c c}
\includegraphics[scale=0.135]{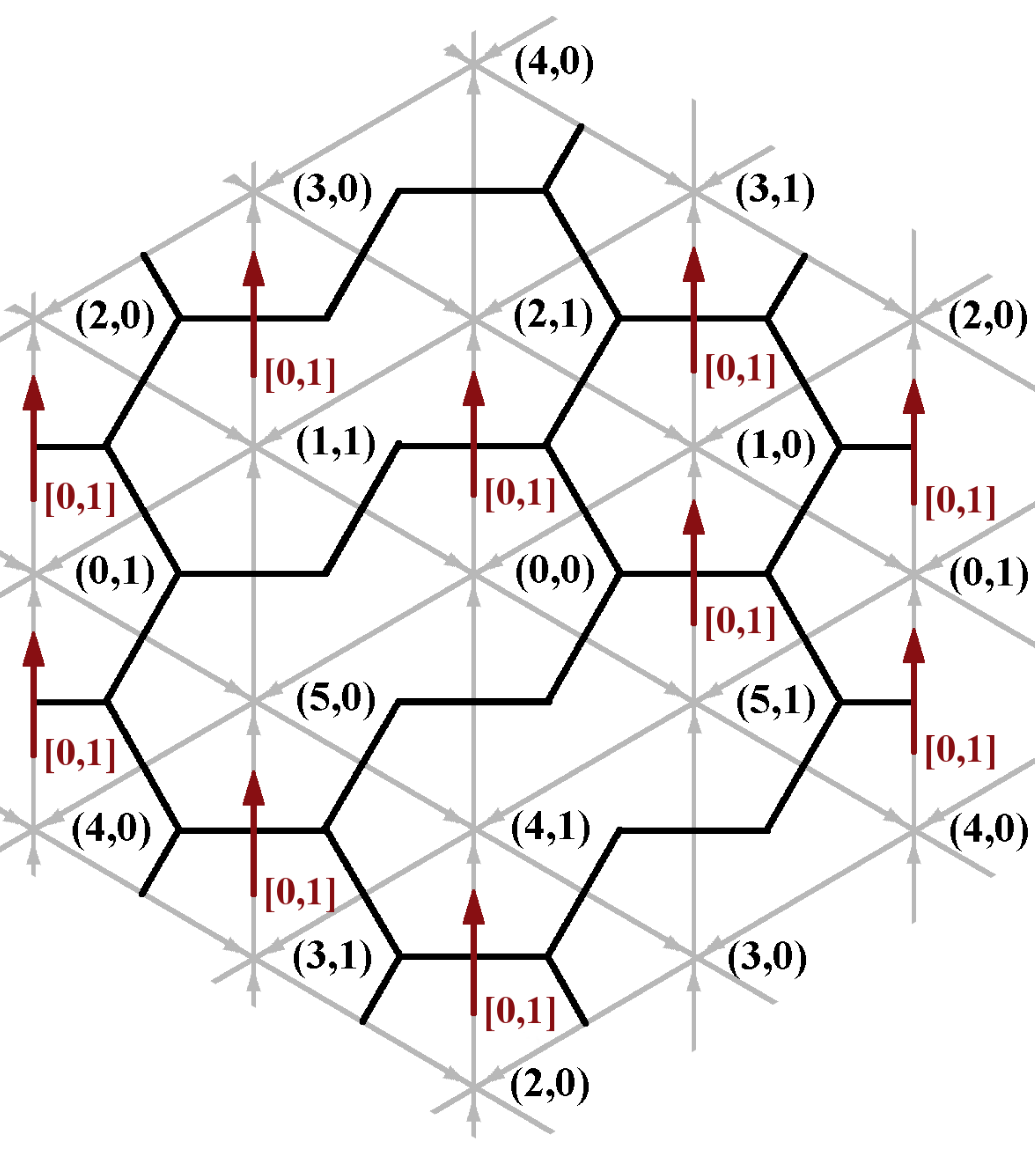} &
\includegraphics[scale=0.135]{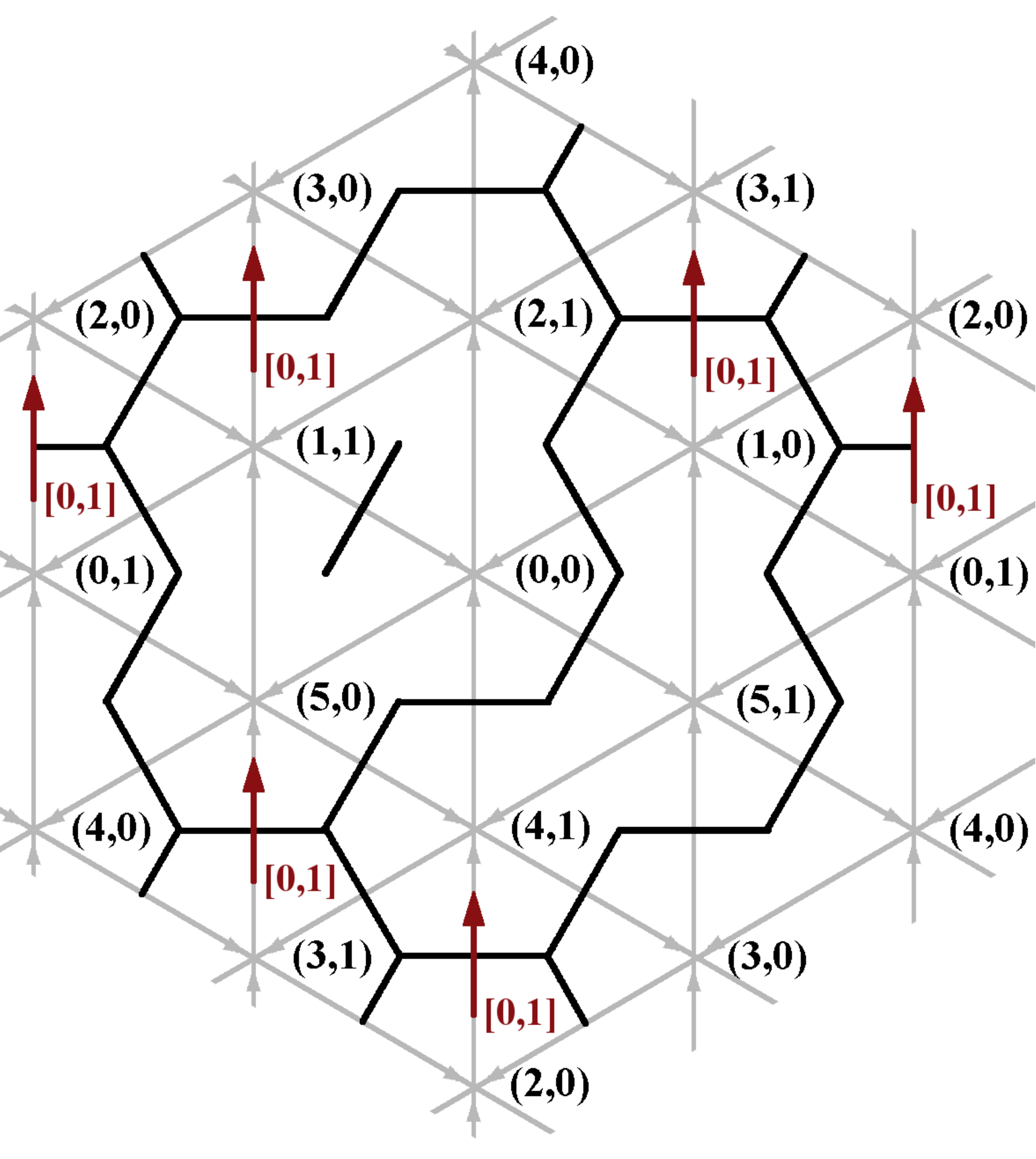} \\
$(S_{3},S_{1,7})$ & $(S_{4},S_{1,7})$ 
\end{tabular}
\begin{tabular}{c c}
\includegraphics[scale=0.135]{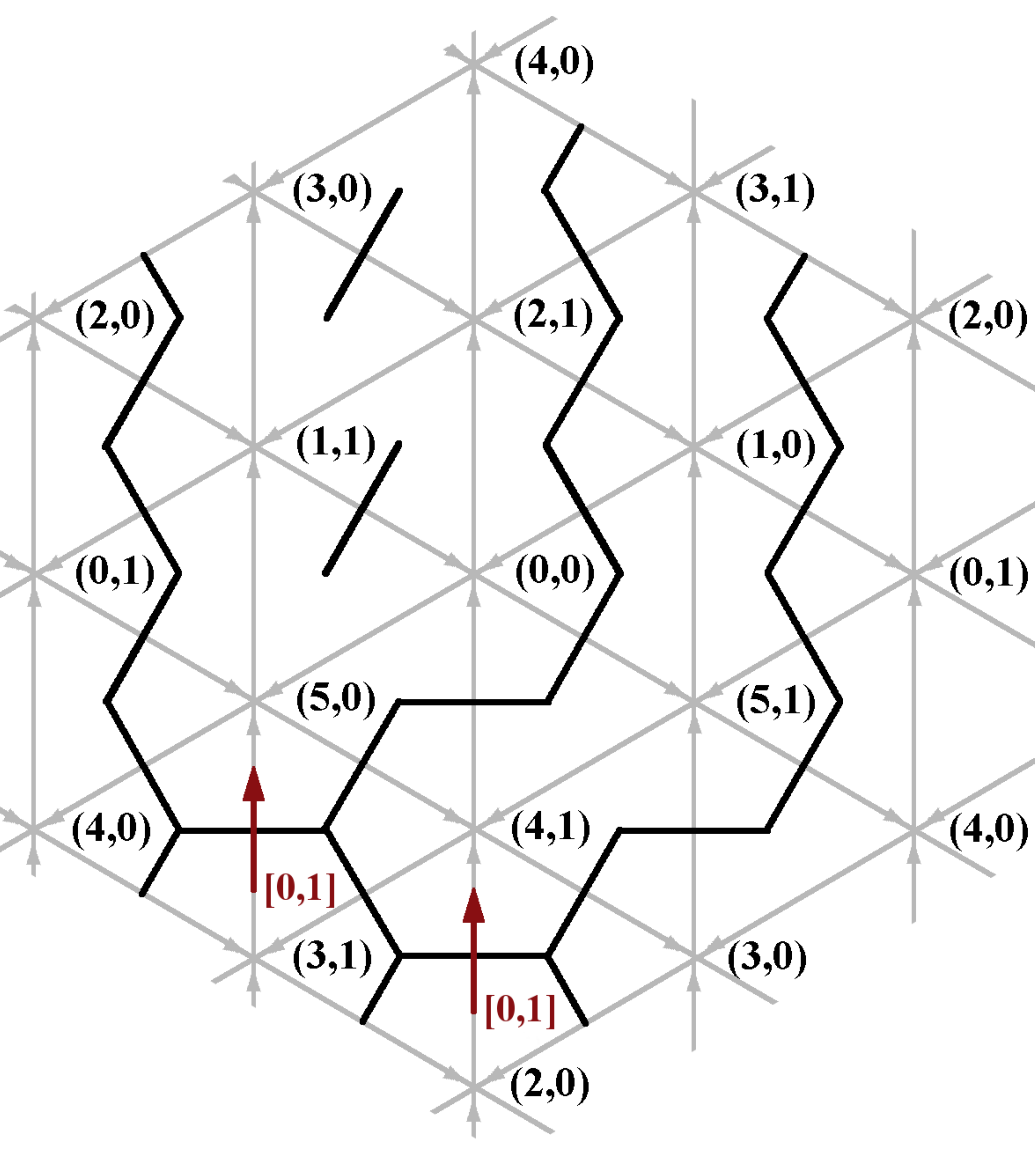} &
\includegraphics[scale=0.135]{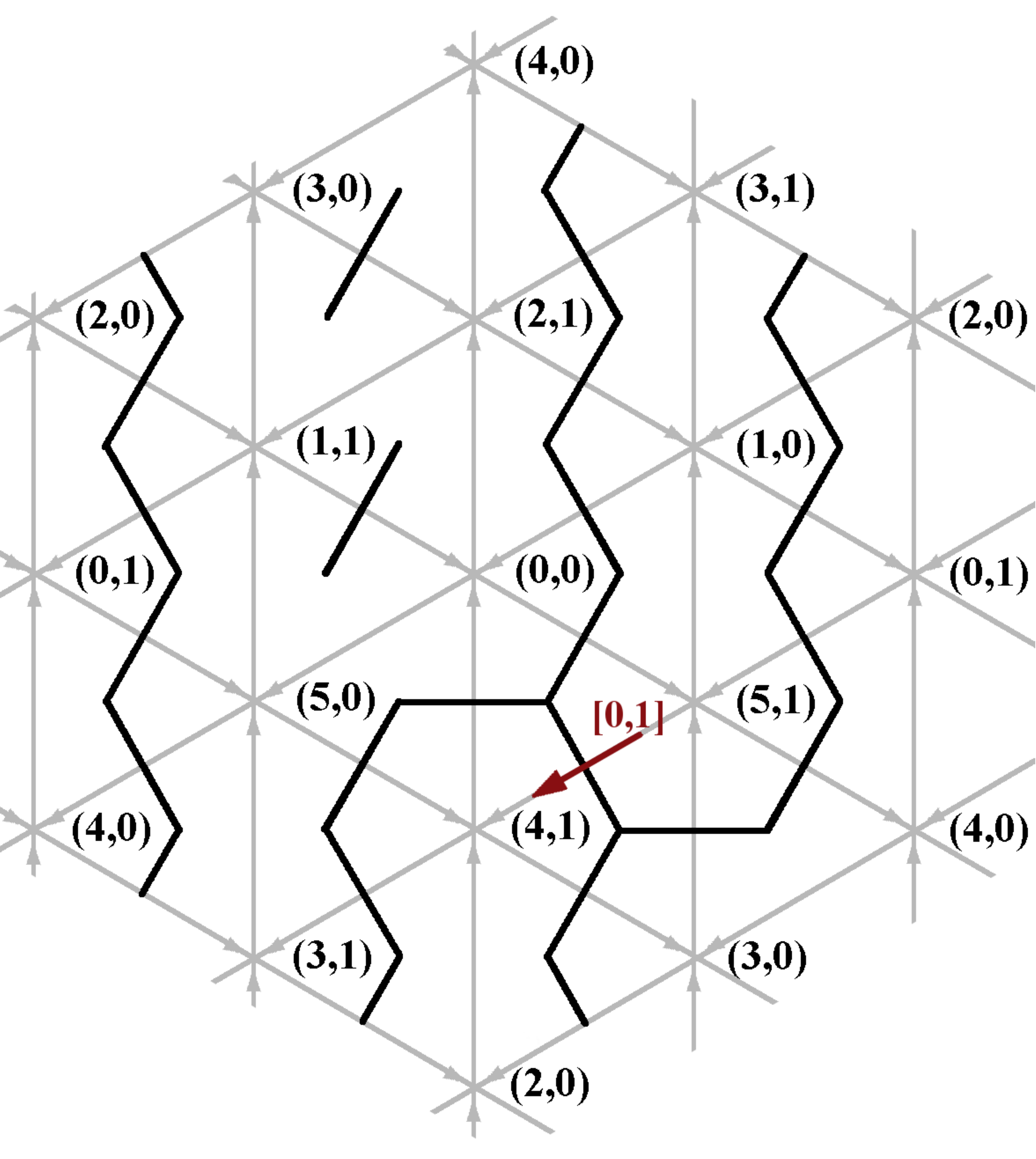} \\
$(S_{5},S_{1,7})$ & $(S_{6},S_{1,7})$ 
\end{tabular}
\begin{tabular}{c c}
\includegraphics[scale=0.135]{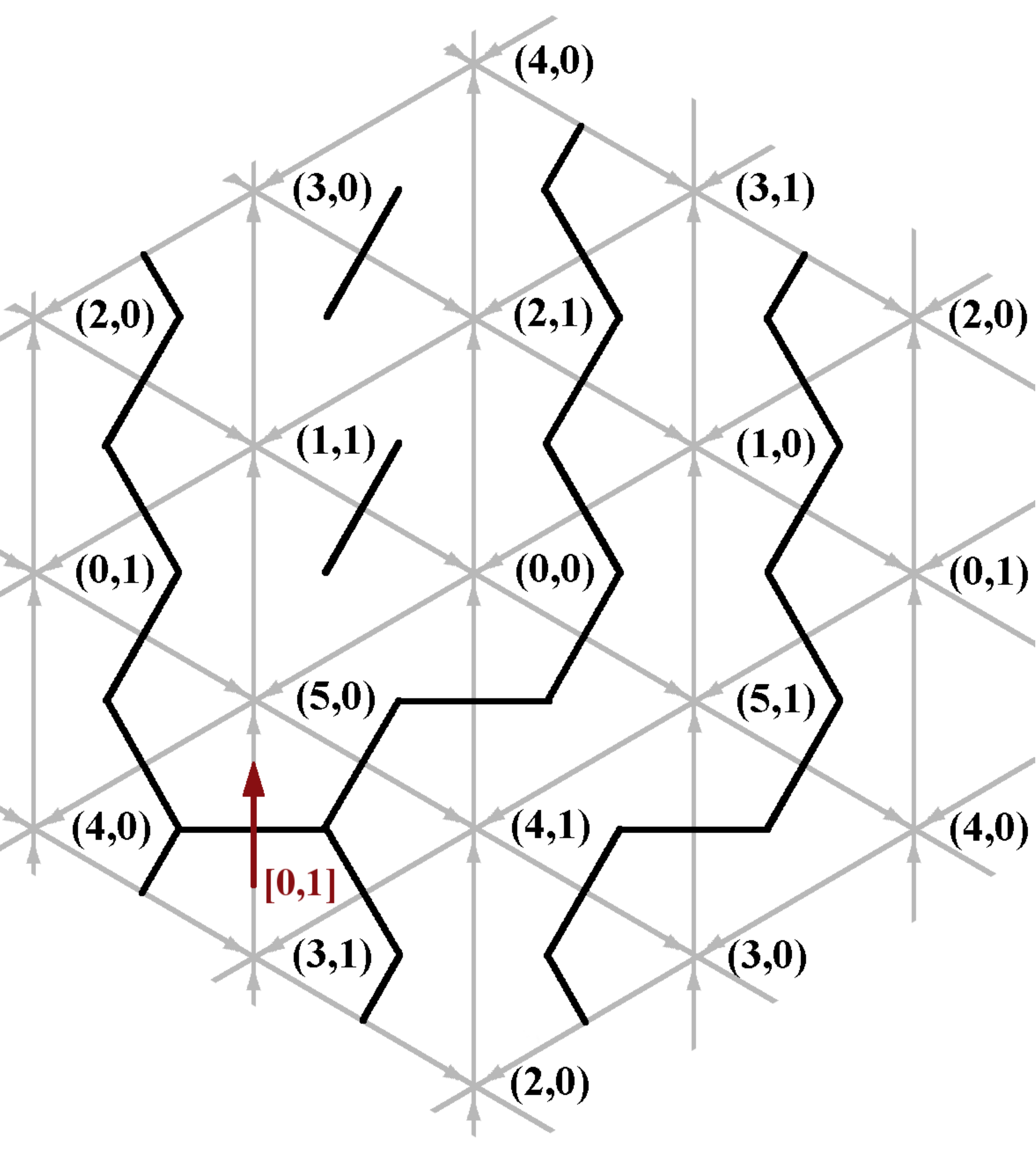} &
\includegraphics[scale=0.135]{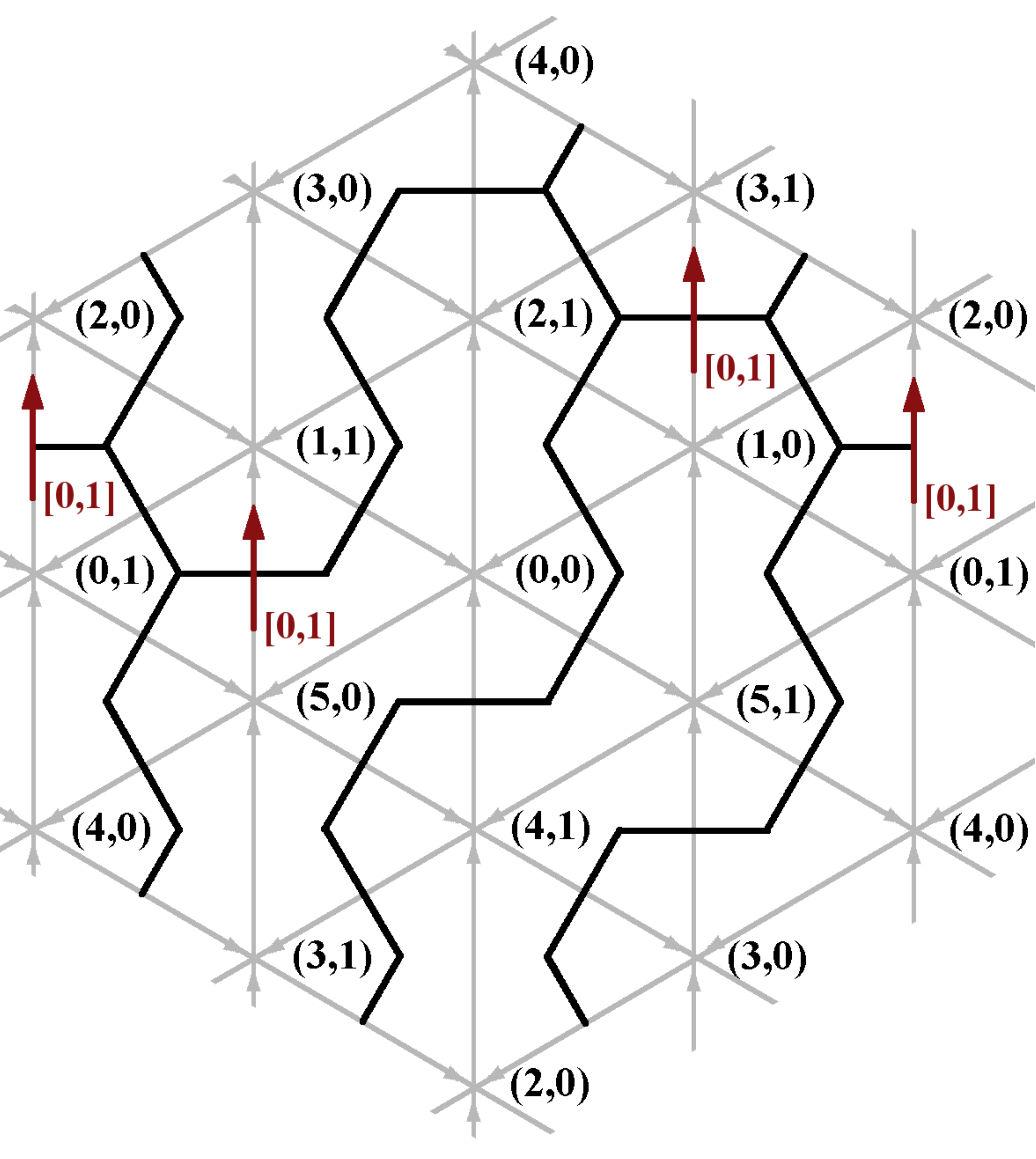} \\
$(S_{9},S_{1,7})$ & $(S_{10},S_{1,7})$ 
\end{tabular}
\end{center}

\bibliography{../../../references}

\providecommand{\bysame}{\leavevmode\hbox to3em{\hrulefill}\thinspace}
\providecommand{\MR}{\relax\ifhmode\unskip\space\fi MR }
\providecommand{\MRhref}[2]{%
  \href{http://www.ams.org/mathscinet-getitem?mr=#1}{#2}
}
\providecommand{\href}[2]{#2}
\begin{thebibliography}{KKMSD73}

\bibitem[BK04]{Kal-Bezr-04}
R.~Bezrukavnikov and D.~Kaledin, \emph{Mc{K}ay equivalence for symplectic
  resolutions of singularities}, Proc. Steklov Inst. Math \textbf{246} (2004),
  13--33, math.AG/0401002.

\bibitem[BKR01]{BKR01}
T.~Bridgeland, A.~King, and M.~Reid, \emph{The {M}c{K}ay correspondence as an
  equivalence of derived categories}, J. Amer. Math. Soc. \textbf{14} (2001),
  535--554, math.AG/9908027.

\bibitem[BM02]{BridgelandMacioca-FourierMukaiTransformsFor_K3_AndEllipticFibra%
tions}
T.~Bridgeland and A.~Macioca, \emph{{F}ourier-{M}ukai transforms for {K}3 and
  elliptic fibrations}, J. Algebraic Geom. \textbf{11} (2002), no.~4, 629--657,
  math.AG/9908022.

\bibitem[BO95]{BonOr95}
A.~Bondal and D.~Orlov, \emph{Semi-orthogonal decompositions for algebraic
  varieties}, preprint math.AG/950612, (1995).

\bibitem[Bri99]{Bridg97}
T.~Bridgeland, \emph{Equivalence of triangulated categories and
  {F}ourier-{M}ukai transforms}, Bull. London Math. Soc. \textbf{31} (1999),
  25--34, math.AG/9809114.

\bibitem[Bri02]{Bridgeland-StabilityConditionsOnTriangulatedCategories}
\bysame, \emph{Stability conditions on triangulated categories}, preprint
  math.AG/0504584, to appear in Annals of Mathematics, (2002).

\bibitem[CI04]{Craw-Ishii-02}
A.~Craw and A.~Ishii, \emph{Flops of ${G}$-$\hilb$ and equivalences of derived
  category by variation of {GIT} quotient}, Duke Math J. \textbf{124} (2004),
  no.~2, 259--307, math.AG/0211360.

\bibitem[CMT05a]{Craw-Maclagan-Thomas-05-I}
A.~Craw, D.~Maclagan, and R.R. Thomas, \emph{Moduli of {M}c{K}ay quiver
  representations {I}: the coherent component}, preprint, (2005).

\bibitem[CMT05b]{Craw-Maclagan-Thomas-05-II}
\bysame, \emph{Moduli of {M}c{K}ay quiver representations {II}: Grobner basis
  techniques}, preprint, (2005).

\bibitem[CR02]{Craw02}
A.~Craw and M.~Reid, \emph{How to calculate ${A}$-$\hilb \mathbb{C}^3$},
  Seminaires et Congres \textbf{6} (2002), 129--154, math.AG/9909085.

\bibitem[Del66]{DeligneCohomologieASupportPropre}
P.~Deligne, \emph{Cohomologie a support propre et construction du foncteur
  $f^!$}, in ``Residues and Duality'', R. Hartshorne, Springer, 1966,
  pp.~404--421.

\bibitem[GD60]{Grothendieck-EGA-I}
A.~Grothendieck and J.~Dieudonn{\'e}, \emph{{\'E}l{\'e}ments de
  g{\'e}om{\'e}trie alg{\'e}brique {I}: Le langage des sch{\'e}mas.},
  Publications math{\'e}matiques de l'I.H.{\'E}.S. \textbf{4} (1960), 5--228.

\bibitem[GD61]{Grothendieck-EGA-III-1}
\bysame, \emph{{\'E}l{\'e}ments de g{\'e}om{\'e}trie alg{\'e}brique {III}:
  {\'E}tude cohomologique des faisceaux coh{\'e}rents. {P}remi{\`e}re partie.},
  Publications math{\'e}matiques de l'I.H.{\'E}.S. \textbf{11} (1961), 5--167.

\bibitem[GD63]{Grothendieck-EGA-III-2}
\bysame, \emph{{\'E}l{\'e}ments de g{\'e}om{\'e}trie alg{\'e}brique {III}:
  {\'E}tude cohomologique des faisceaux coh{\'e}rents. {S}econde partie.},
  Publications math{\'e}matiques de l'I.H.{\'E}.S. \textbf{17} (1963), 5--91.

\bibitem[GD66]{Grothendieck-EGA-IV-3}
\bysame, \emph{{\'E}l{\'e}ments de g{\'e}om{\'e}trie alg{\'e}brique {IV}:
  {\'E}tude locale des sch{\'e}mas et des morphismes de sch{\'e}mas.
  {T}roisi{\`e}me partie.}, Publications math{\'e}matiques de l'I.H.{\'E}.S.
  \textbf{28} (1966), 5--255.

\bibitem[GM03]{GelfandManin-MethodsOfHomologicalAlgebra}
S.I. Gelfand and Yu.~I. Manin, \emph{Methods of homological algebra}, Springer,
  2003.

\bibitem[GSV83]{GsV-ConstructionGeometriqueDeLaCorrespondanceDeMcKay}
G.~Gonzales-Sprinberg and J.-L. Verdier, \emph{Construction g{\'e}om{\'e}trique
  de la correspondance de {M}c{K}ay}, Ann. sci. ENS \textbf{16} (1983),
  409--449.

\bibitem[Har66]{Hartshorne-Residues-and-Duality}
R.~Hartshorne, \emph{Residues and duality}, Springer-Verlag, 1966.

\bibitem[Huy06]{Huybrechts-FourierMukaiTransformsInAlgebraicGeometry}
D.~Huybrechts, \emph{{F}ourier-{M}ukai transforms in algebraic geometry},
  Oxford University Press, 2006.

\bibitem[IN00]{ItoNakajima98}
Y.~Ito and H.~Nakajima, \emph{Mc{K}ay correspondence and {H}ilbert schemes in
  dimension three}, Topology \textbf{39} (2000), no.~6, 1155--1191,
  math.AG/9803120.

\bibitem[Kaw05]{Kawamata-LogCrepantBirationalMapsAndDerivedCategories}
Y.~Kawamata, \emph{Log crepant birational maps and derived categories}, J.
  Math. Sci. Univ. Tokyo \textbf{12} (2005), no.~2, 211--231, math.AG/0311139.

\bibitem[Kin94]{King94}
A.~King, \emph{Moduli of representations of finite-dimensional algebras},
  Quart. J. Math. Oxford \textbf{45} (1994), 515--530.

\bibitem[KKMSD73]{KKMSD-ToroidalEmbeddingsI}
G.~Kempf, F.~Knudsen, D.~Mumford, and B.~Saint-Donat, \emph{Toroidal embeddings
  {I}}, Springer-Verlag, 1973.

\bibitem[Kol89]{Kollar-Flops}
J.~Koll\'ar, \emph{Flops}, Nagoya Math J. \textbf{113} (1989), 15--36.

\bibitem[Kuz05]{Kuznetsov-HomologicalProjectiveDuality}
A.~Kuznetsov, \emph{Homological projective duality}, preprint math.AG/0507292,
  (2005).

\bibitem[KV98]{KapranovVasserot-KleinianSingularitiesDerivedCategoriesAndHallA%
lgebras}
M.~Kapranov and E.~Vasserot, \emph{Kleinian singularities, derived categories
  and hall algebras}, preprint math.AG/9812016, (1998).

\bibitem[Log03]{Logvinenko-Families-of-G-constellations-over-resolutions-of-qu%
otient-singularities}
T.~Logvinenko, \emph{Families of {$G$}-constellations over resolutions of
  quotient singularities}, preprint math.AG/0305194, (2003).

\bibitem[Log04]{Logvinenko-thesis}
\bysame, \emph{Families of ${G}$-{C}onstellations parametrised by resolutions
  of quotient singularities}, Ph.D. thesis, University of {B}ath, 2004.

\bibitem[Log06]{Logvinenko-Natural-G-Constellation-Families}
\bysame, \emph{Natural {$G$}-constellation families}, preprint math.AG/0601014,
  (2006).

\bibitem[Mat86]{Mats86}
H.~Matsumura, \emph{Commutative ring theory}, Cambridge University Press, 1986.

\bibitem[McK80]{McKay-GraphsSingularitiesAndFiniteGroups}
J.~McKay, \emph{Graphs, singularities and finite groups}, Proc. Symp. Pure
  Math. \textbf{37} (1980), 183--186.

\bibitem[Muk81]{Mukai-DualityBetweenDXandDX^WithItsApplicationToPicardSheaves}
S.~Mukai, \emph{Duality between {$D(X)$} and {$D(\hat{X})$} with its
  application to {P}icard sheaves}, Nagoya Math J \textbf{81} (1981), 153--175.

\bibitem[Nag62]{Nagata-ImbeddingOfAnAbstractVarietyInACompleteVariety}
M.~Nagata, \emph{Imbedding of an abstract variety in a complete variety}, J.
  Math. Kyoto Uni. \textbf{2} (1962), no.~1.

\bibitem[Orl97]{Orlov-EquivalencesOfDerivedCategoriesAndK3Surfaces}
D.~Orlov, \emph{Equivalences of derived categories and {K}3 surfaces}, J. Math.
  Sci. (NY) \textbf{84} (1997), no.~5, 1361--1381, math.AG/9606006.

\bibitem[Rei87]{YPG87}
M.~Reid, \emph{Young person's guide to canonical singularities}, Proc. of
  Symposia in Pure Math. \textbf{46} (1987), 345--414.

\bibitem[Rei97]{Kinosaki-97}
\bysame, \emph{{M}c{K}ay correspondence}, preprint math.AG/9702016, (1997).

\bibitem[Rob98]{Roberts-MultiplicitiesAndChernClassesInLocalAlgebra}
P.~C. Roberts, \emph{Multiplicities and {C}hern classes in local algebra},
  Cambridge University Press, 1998.

\bibitem[Ser00]{Serre-LocalAlgebra}
J.P. Serre, \emph{Local algebra}, Springer, 2000.

\end{thebibliography}
\bibliographystyle{amsalpha}
\end{document}